# PATHWISE UNIQUENESS FOR A DEGENERATE STOCHASTIC DIFFERENTIAL EQUATION[1]

By Richard F. Bass, Krzysztof Burdzy and Zhen-Qing Chen

*University of Connecticut, University of Washington
and University of Washington*

We introduce a new method of proving pathwise uniqueness, and we apply it to the degenerate stochastic differential equation
$$dX_t = |X_t|^\alpha \, dW_t,$$
where $W_t$ is a one-dimensional Brownian motion and $\alpha \in (0, 1/2)$. Weak uniqueness does not hold for the solution to this equation. If one restricts attention, however, to those solutions that spend zero time at 0, then pathwise uniqueness does hold and a strong solution exists. We also consider a class of stochastic differential equations with reflection.

**1. Introduction.** In this paper we introduce a new method of proving pathwise uniqueness for certain stochastic differential equations. The technique uses ideas from excursion theory. We apply this method to the degenerate stochastic differential equation

(1.1) $$dX_t = |X_t|^\alpha \, dW_t,$$

where $W_t$ is a one-dimensional Brownian motion. When $\alpha \in [1/2, 1]$, the classical theorem of Yamada–Watanabe [12] says that pathwise uniqueness holds for (1.1). Moreover, this is sharp: it is well known that pathwise uniqueness does not hold for (1.1) when $\alpha \in (0, 1/2)$. In fact, even weak uniqueness (i.e., uniqueness in law) does not hold when $\alpha \in (0, 1/2)$. When $x_0 = 0$, one solution is the identically zero one, while a nonzero solution can be constructed by time changing a Brownian motion.

This, however, is not all that can be said about uniqueness for (1.1). One of the main points of this paper is that the only reason pathwise uniqueness

Received January 2006; revised November 2006.
[1]Supported in part by NSF Grants DMS-02-44737 and DMS-03-03310.

fails in (1.1) is that weak uniqueness fails. It was shown by Engelbert and Hess [6] and Engelbert and Schmidt [7] that for every $x_0 \in \mathbb{R}$, there is a weak solution to (1.1) that spends zero time at 0 and the law of such a solution is unique. In this paper we show that there is pathwise uniqueness among those solutions to (1.2) that spend zero time at 0 and a strong solution exists.

Before we give rigorous statements of our main results, we recall some definitions.

DEFINITION 1.1.   (i) Given a Brownian motion $W$ on a probability space, a strong solution to the stochastic differential equation

$$(1.2) \qquad X_t = x + \int_0^t |X_s|^\alpha \, dW_s$$

that spends zero time at 0 is a continuous process $X = \{X_t, t \geq 0\}$ that is adapted to the filtration generated by $W$, solves (1.2), and satisfies

$$(1.3) \qquad \int_0^\infty \mathbf{1}_{\{0\}}(X_s) \, ds = 0 \qquad \text{a.s.}$$

(ii) A weak solution of (1.2) is a couple $(X, W)$ on a filtered probability space $(\Omega, \mathcal{F}, \{\mathcal{F}_t\}_{t \geq 0}, \mathbb{P})$ such that $X_t$ is adapted to $\mathcal{F}_t$, $W_t$ is an $\{\mathcal{F}_t\}_{t \geq 0}$-Brownian motion (i.e., $W_t$ is $\mathcal{F}_t$-measurable and for $t > s$, $W_t - W_s$ is independent of $\mathcal{F}_s$ and has a normal distribution with zero mean and variance $t - s$), and $(X, W)$ satisfies (1.2).

(iii) We say weak uniqueness holds for (1.2) among solutions that spend zero time at 0 if whenever $(X, W)$, $(\widetilde{X}, \widetilde{W})$ are two weak solutions of (1.2) satisfying the condition (1.3), then the process $X = \{X_t, t \geq 0\}$ has the same law as the process $\widetilde{X} = \{\widetilde{X}_t, t \geq 0\}$.

(iv) Pathwise uniqueness is said to hold for (1.2) among solutions that spend zero time at 0 if whenever $(X, W)$, $(\widetilde{X}, W)$ are two weak solutions of (1.2) satisfying (1.3) with a common Brownian motion $W$ (relative to possibly different filtrations) on a common probability space and with common initial value, then $\mathbb{P}(X_t = \widetilde{X}_t \text{ for all } t \geq 0) = 1$.

(v) Strong uniqueness is said to hold for (1.2) among solutions that spend zero time at 0 if whenever $(X, W)$, $(\widetilde{X}, W)$ are two weak solutions of (1.2) satisfying (1.3) with a common Brownian motion $W$ on a common probability filtered space and with common initial value, then $\mathbb{P}(X_t = \widetilde{X}_t \text{ for all } t \geq 0) = 1$.

It is clear that pathwise uniqueness implies strong uniqueness. We warn the reader that what we call "strong uniqueness" is sometimes called "pathwise uniqueness," for example, in [11], Definition IX.1.3. We follow [1] in distinguishing between strong uniqueness and pathwise uniqueness. We note that strong uniqueness implies weak uniqueness, by the same argument as in [11], Theorem X1.7(i).

Our main theorem is the following.

THEOREM 1.2. *Suppose $\alpha \in (0, \frac{1}{2})$ and $x \in \mathbb{R}$. Then pathwise uniqueness holds for solutions of (1.2) that spend zero time at 0. Moreover, a strong solution to (1.2) which spends zero time at 0 exists.*

In the above theorem we have both $X$ and $X'$ satisfying (1.2) with respect to the same Brownian motion, but we allow the possibility that there are two different filtrations $\{\mathcal{F}_t\}$ and $\{\mathcal{F}'_t\}$; the process $W$ must be a Brownian motion with respect to both filtrations. The fact that we allow the filtrations to be different does not increase the generality of the theorem in a substantial way—see the last paragraph of the proof of Theorem 1.2.

The lack of weak uniqueness is not the only reason pathwise uniqueness can fail. Barlow [1] showed that for any $\beta < 1/2$ one can construct a bounded Hölder continuous function $\sigma_\beta$ of order $\beta$ that is bounded below such that pathwise uniqueness fails for $dX_t = \sigma_\beta(X_t)\, dW_t$. Weak uniqueness does hold in this case. For other positive results on pathwise uniqueness, see Nakao [10] and Le Gall [8].

Our method of proof of Theorem 1.2 is new and substantially different from any of the existing methods of proving pathwise uniqueness. Some of these previous methods include an appropriate use of Itô's formula, a study of local times, looking at the maximum or minimum of two solutions and constructing a strong solution. We were unable to successfully adapt any of these methods to the study of (1.1).

At the basis of our new method are ideas from excursion theory. We first show that if $X$ and $Y$ are two solutions, $-X_0 \le Y_0 \le X_0$, and $X$ is conditioned to hit the level 1 before hitting the level 0, then when $X$ hits the level 1, the process $Y$ will also be close to the level 1 with high probability, provided $X_0$ is small enough. We refer to this as the "chasing phenomenon." We then use this to show that for every $\delta > 0$, with probability one, the processes $X$ and $Y$ have to agree on every excursion of $M := |X| \vee |Y|$ away from zero that reaches level $\delta$, which establishes the pathwise uniqueness.

Interestingly, the one-sided problem, that is, pathwise uniqueness for stochastic differential equations with reflection, is much easier. Consider the equation

$$(1.4) \qquad X_t = x + \int_0^t a(X_s)\, dW_s + \int_0^t b(X_s)\, ds + L_t,$$

where $L_t$ is a nondecreasing continuous process that increases only when $X$ is at 0, $X_t$ is never negative, and

$$(1.5) \qquad \int_0^\infty \mathbf{1}_{\{0\}}(X_s)\, ds = 0 \qquad \text{a.s.}$$

We define pathwise uniqueness, strong solution, and weak solution for (1.4) analogously to Definition 1.1. When $a(x) = |x|^\alpha$ with $\alpha \in (0, 1/2)$ and $b(x) \equiv 0$, pathwise uniqueness was proved in [5]. We give a theorem for solutions of (1.4) that greatly generalizes the result of [5], with a much simpler proof.

THEOREM 1.3. *Suppose that $b$ is a bounded measurable function on $\mathbb{R}$. Suppose that the function $a:[0,\infty) \to [0,\infty)$ is bounded, $a^{-2}$ is locally integrable on $\mathbb{R}$ and satisfies either:*

(a) *for every closed subinterval $I$ of $(0,\infty)$ there exists a continuous increasing function $\rho_I:[0,\infty) \to [0,\infty)$ with $\rho_I(0) = 0$ such that*

$$(1.6) \qquad |a(x) - a(y)| \le \rho_I(|x-y|), \qquad x,y \in I,$$

*and*

$$(1.7) \qquad \int_{0+} \frac{1}{\rho_I(h)^2}\, dh = \infty,$$

*or*

(b) *on each closed subinterval $I$ of $(0,\infty)$ the coefficient $a$ is bounded below by a positive constant and is of finite quadratic variation.*

*Then pathwise uniqueness holds for solutions of (1.4) that spend zero time at 0. Moreover there is a strong solution to (1.4) that spends zero time at 0.*

As an immediate application of Theorem 1.3, we have the following.

COROLLARY 1.4. *Suppose that $b$ is an odd bounded measurable function on $\mathbb{R}$. Suppose that $a$ is an odd bounded measurable function on $\mathbb{R}$ with $a^{-2}$ locally integrable on $\mathbb{R}$ and satisfies either condition* (a) *or* (b) *in Theorem 1.3. Then for any two weak solutions $(X, W)$ and $(\widetilde{X}, \widetilde{W})$ to*

$$(1.8) \qquad X_t = x + \int_0^t a(X_s)\, dW_s + \int_0^t b(X_s)\, ds$$

*with*

$$(1.9) \qquad \int_0^\infty \mathbf{1}_{\{0\}}(X_s)\, ds = 0 \qquad a.s.,$$

*with a common Brownian motion $W$ (relative to possibly different filtrations) on a common probability space and with common initial value, we have*

$$\mathbb{P}(|X_t| = |\widetilde{X}_t| \text{ for all } t \ge 0) = 1.$$

The above corollary extends the main result (Theorem 1) of [9].

REMARK 1.5. We do not fully understand why the proof of Theorem 1.3 is so much easier than that of Theorem 1.2. The proof of the one-sided version of Theorem 1.2 given in [5] proceeds by constructing a strong solution to (1.2). If one tries that in the two-sided context, one gets the difference of two terms each tending to infinity, and one is not able to prove convergence. The proof of Theorem 1.3 given here also does not extend to the two-sided context. In addition, it is easy to see that one cannot derive Theorem 1.2 just by applying Theorem 1.3 to $|X_t|, |Y_t|$. See Remark 5.1 in Section 5 for more comments.

Let us mention two open problems which we think are quite interesting.

PROBLEM 1. Consider the equation

$$dX_t = a(X_t)\,dW_t$$

and let $Z_a := \{x : a(x) = 0\}$. Suppose $a$ is smooth on every closed interval contained in $Z_a^c$, $a^{-2}$ is locally integrable, and $a$ is bounded. Is there pathwise uniqueness among those solutions that spend zero time in the set $Z_a$? The smoothness is needed to rule out counter examples such as those of Barlow [1]. The local integrability is necessary for a weak solution that spends zero time at $Z_a$ to be unique; see [7].

PROBLEM 2. For each $\lambda \in [0,\infty)$ there is a strong Markov process $X$ satisfying (1.2) and associated to the speed measure

$$m(dx) = |x|^{-2\alpha}\,dx + \lambda\delta_0(dx),$$

where $\delta_0$ is point mass at 0. The value $\lambda$ measures how "sticky" the diffusion is at 0. Theorem 1.2 covers the case $\lambda = 0$. What can one say for other values of $\lambda$? Uniqueness in law holds for each value of $\lambda$. When $\lambda \neq 0$, what additional condition or conditions must one impose on solutions to $dX_t = |X_t|^\alpha\,dW_t$ so that the solutions have the speed measure given above? Does pathwise uniqueness hold in this situation?

In the next section we discuss some preliminaries. Section 3 discusses the chasing phenomenon, while the proof of Theorem 1.2 is given in Section 4. Theorem 1.3 is proved in Section 5.

**2. Preliminaries.** Suppose that $X$ and $Y$ are two weak solutions to (1.1) and (1.3) driven by the same Brownian motion $W$, starting from $W_0 = w$, $X_0 = x$ and $Y_0 = y$, and defined on some probability space $(\Omega, \mathbb{P})$. Let $(\Omega^C, \mathcal{F}^C, \{\mathcal{F}_t^{0,C}\}_{t\geq 0})$ be the canonical probability space, that is, $\Omega^C$ is the collection of continuous functions from $[0,\infty)$ to $\mathbb{R}^3$. For $\omega \in \Omega^C$, we write $\omega = (\omega_1, \omega_2, \omega_3)$, $W_t^C(\omega) = \omega_1(t)$, $X_t^C(\omega) = \omega_2(t)$ and $Y_t^C(\omega) = \omega_3(t)$. The $\sigma$-field $\mathcal{F}^C$ is generated by the cylindrical sets, and $\{\mathcal{F}_t^{0,C}\}_{t \geq 0}$ is the natural filtration generated by $(W^C, X^C, Y^C)$. We now define $\mathbb{P}^{w,x,y}$ to be the probability on the space $(\Omega^C, \mathcal{F}^C)$ such that for every $n \geq 1$, every Borel measurable subset $A$ of $\mathbb{R}^{3n}$, and every choice $t_1 \leq t_2 \leq \cdots \leq t_n$ we have

$$\begin{aligned}
\mathbb{P}^{w,x,y}&(((W_{t_1}^C, X_{t_1}^C, Y_{t_1}^C), \ldots, (W_{t_n}^C, X_{t_n}^C, Y_{t_n}^C)) \in A) \\
&= \mathbb{P}(((W_{t_1}, X_{t_1}, Y_{t_1}), \ldots, (W_{t_n}, X_{t_n}, Y_{t_n})) \in A).
\end{aligned} \tag{2.1}$$

At this point, we cannot assume that for each fixed $(w, x, y)$, the joint law of $(W, X, Y)$ for two weak solutions $X$ and $Y$ to (1.1) and (1.3) driven by the same Brownian motion $W$ with $(W_0, X_0, Y_0) = (w, x, y)$ is unique—this is what we will prove in this paper. Hence, we let $\mathcal{P}(w, x, y)$ denote the collection of all measures $\mathbb{P}^{w,x,y}$ on $(\Omega^C, \mathcal{F}^C)$ obtained by the above recipe. Each triple $(W, X, Y)$ of weak solutions $X$ and $Y$ to (1.1) and (1.3), driven by the same Brownian motion $W$, and starting from $W_0 = w$, $X_0 = x$ and $Y_0 = y$ will give rise to an element of $\mathcal{P}(w, x, y)$. For every measure $\mathbb{P}^{w,x,y} \in \mathcal{P}(w, x, y)$, it is easy to construct distinct triples $(W, X, Y)$ and $(W', X', Y')$ corresponding to $\mathbb{P}^{w,x,y}$, for example, by defining the processes $(W, X, Y)$ and $(W', X', Y')$ on different probability spaces. Whenever we make an assertion about $\mathbb{P}^{w,x,y}$, it should be understood that it holds for all $\mathbb{P}^{w,x,y} \in \mathcal{P}(w, x, y)$.

Most of the time, the value of the index $w$ in $\mathbb{P}^{w,x,y}$ will be irrelevant. Hence we will write $\mathbb{P}^{x,y}$ instead of $\mathbb{P}^{w,x,y}$. Any assertion made about $\mathbb{P}^{x,y}$ should be understood as an assertion that applies to all $\mathbb{P}^{w,x,y} \in \mathcal{P}(w, x, y)$, for all values of $w$. Thus we will abbreviate our notation by referring to $\mathbb{P}^{x,y} \in \mathcal{P}(x, y)$ rather than $\mathbb{P}^{w,x,y} \in \mathcal{P}(w, x, y)$.

Note that under $\mathbb{P}^{x,y} \in \mathcal{P}(x, y)$, $X^C$ and $W^C$ satisfy (1.2)–(1.3) because the stochastic integral can be defined as an almost sure limit along a sequence of discrete approximations, and the finite-dimensional distributions for $(W^C, X^C)$ are the same as those for $(W, X)$, by (2.1). As we mentioned above, [6, 7] prove that for every $x$ there exists a weak solution to (1.2) and (1.3), and that weak uniqueness holds for these solutions. Hence, if $A \in \sigma(X_t^C, t \geq 0)$, then for any $x, y_1$ and $y_2$, we have $\mathbb{P}^{x,y_1}(A) = \mathbb{P}^{x,y_2}(A)$. Therefore, for events $A \in \sigma(X_t^C, t \geq 0)$, we will write $\mathbb{P}^x(A)$ to indicate that $\mathbb{P}^x(A)$ is the common value of $\mathbb{P}^{x,y}(A)$ for all $y$.

Our goal in this paper is to show that if $X$ and $Y$ are weak solutions to (1.2) and (1.3) with the same initial value $x$ and driven by the same Brownian motion $W$, then almost surely $X_t = Y_t$ for all $t \geq 0$. In view of (2.1) it suffices to prove that $\mathbb{P}^{x,x}(X_t^C = Y_t^C \text{ for all } t) = 1$. We can thus restrict our attention to the canonical probability space and $(W^C, X^C, Y^C)$. We do so henceforth, and we drop the superscript "$C$" from now on.

We define the "minimal augmented filtration" $\{\mathcal{F}_t\}_{t \geq 0}$ on $\Omega$ as follows. For each $\mathbb{P}^{x,y} \in \mathcal{P}(x, y)$ we add all $\mathbb{P}^{x,y}$-null sets to each $\mathcal{F}_{t+}^0$. If we denote the $\sigma$-field so formed by $\mathcal{F}_t(\mathbb{P}^{x,y})$, we then form the filtration where

$$\mathcal{F}_t = \bigcap_{\mathbb{P}^{x,y} \in \mathcal{P}(x,y)} \mathcal{F}_t(\mathbb{P}^{x,y}).$$

The minimal augmented filtration $\{\mathcal{F}_t\}_{t \geq 0}$ that we just defined is right continuous. We define $\mathcal{F} = \mathcal{F}_\infty := \sigma(\mathcal{F}_t, t \geq 0)$.

For $a \in \mathbb{R}$ we will write $T_a^X = \inf\{t > 0 : X_t = a\}$. Similar notation will be used for hitting times of other processes. When there is no confusion possible, we will write $T_a$ for $T_a^X$. We will often use the following stopping time: $T = T_0^X \wedge T_1^X$.

We will sometimes look at $\mathbb{P}^x(A | T_1^X < T_0^X)$, where $A \in \mathcal{F}_{T_1^X \wedge T_0^X}$. We will explain now how this probability can be represented using Doob's $h$-transform. Let $\Delta$ denote a cemetery state added to the state space of $X$. For an open interval $D \subset [0, \infty)$, let $X^D$ denote the process $X$ killed upon leaving $D$, that is, $X_t^D = X_t$ for $t < T_{D^c}$, and $X_t^D = \Delta$ for $t \geq T_{D^c}$. Since $X$ is on natural scale, $\mathbb{P}^x(T_1^X < T_0^X) = x$ for $x \in (0,1)$, and $h(x) = x$ is a harmonic function for $X^{(0,1)}$. We define the conditional law $\mathbb{Q}_1^x$ of $X^{(0,1)}$ starting from $x \in (0,1)$ given the event $\{T_1^X < T_0^X\}$ by

$$\mathbb{Q}_1^x(A) = \frac{1}{h(x)} \mathbb{E}^x[\mathbf{1}_A h(X_{T_0^X \wedge T_1^X})],$$

for $A \in \sigma(X_t, t \geq 0)$ such that $A \in \mathcal{F}_{T_1^X \wedge T_0^X}$. In fact, $\mathbb{Q}_1^x(A)$ is well defined for any $A \in \sigma(X_t, t \geq 0)$ by the above formula and later in this paper we sometimes do take such an extension. We use $\mathbb{Q}_0^x$ to denote the law of $h$-transformed process $X^{(0,1)}$ starting from $x$ with $h(y) = 1 - y$; this corresponds to $X$ starting from $x$ conditioned to hit $0$ before $1$. We use $\mathbb{Q}_\infty^x$ for the law of $X^{(0,\infty)}$ $h$-path transformed by the function $h(y) = y$.

According to our conventions, $\mathbb{Q}_1^x(A)$ is a special case of $\mathbb{Q}_1^{x,y}(A)$ when $A \in \sigma(X_t, t \geq 0)$, where $\mathbb{Q}_1^{x,y}(A)$ is defined for all $A \in \mathcal{F}$ by

(2.2) $$\mathbb{Q}_1^{x,y}(A) = \frac{1}{h(x)} \mathbb{E}^{x,y}[\mathbf{1}_A h(X_{T_0^X \wedge T_1^X})].$$

We let $\mathcal{Q}_1(x,y)$ be the collection of all such $\mathbb{Q}_1^{x,y}$ when $\mathbb{P}^{x,y} \in \mathcal{P}(x,y)$. We similarly define $\mathbb{Q}_0^{x,y}, \mathbb{Q}_\infty^{x,y}, \mathcal{Q}_0(x,y)$ and $\mathcal{Q}_\infty(x,y)$.

Recall the usual shift operator notation $\theta_t$. For example, we write $u + T_a^X \circ \theta_u$ for the stopping time $\inf\{t > u : X_t = a\}$.

REMARK 2.1. As indicated previously, we cannot assume that the pair $(X, Y)$ is strong Markov. However, there is a substitute for the strong Markov property that is almost as useful. For a probability $\overline{\mathbb{P}} \in \mathcal{P}(x,y)$, let $S$ be a finite stopping time with respect to the filtration $\{\mathcal{F}_t, t \geq 0\}$ and define $\overline{\mathbb{P}}_S(A) = \overline{\mathbb{P}}(A \circ \theta_S)$ for $A \in \mathcal{F}_\infty$. Let $\mathbb{O}_S(\cdot, \cdot)$ be a regular conditional probability for $\overline{\mathbb{P}}_S$ given $\mathcal{F}_S$. Then if $X$ and $Y$ are two solutions to (1.1) and (1.3), starting from $X_0$ and $Y_0$, respectively, and driven by the same Brownian motion, under $\mathbb{O}_S$ the processes $(X_{S+t}, Y_{S+t})$ are again solutions to (1.1) and (1.3) driven by the same Brownian motion that spend zero time at 0, started at $(X_S, Y_S)$. In other words, $\mathbb{O}_S \in \mathcal{P}(X_S, Y_S)$ a.s., if $\overline{\mathbb{P}} \in \mathcal{P}(x, y)$. The proof of this is the same as the proof of Proposition VI.2.1 in [3], except for showing that zero time is spent at 0. This last fact follows easily because

$$\mathbb{E}\left[\mathbb{E}_{\mathbb{O}_S}\left[\int_0^\infty \mathbf{1}_{\{0\}}(X_s)\,ds\right]\right] = \mathbb{E}\left[\int_S^\infty \mathbf{1}_{\{0\}}(X_s)\,ds\right] = 0,$$

hence $\mathbb{O}_S(\int_0^\infty \mathbf{1}_{\{0\}}(X_s)\,ds \neq 0)$ is zero for almost every $\omega$; the same argument applies to $Y$. We refer to this as the pseudo-strong Markov property. See also [4] for examples as to how the pseudo-strong Markov property is used. We will give the full argument in our first nontrivial use of this property below (see the proof of Lemma 3.4), but in other usages leave the details to the reader.

REMARK 2.2. Recall that $T = T_0^X \wedge T_1^X$. If $x \in (0, 1)$, the measure $\mathbb{P}^x$ is the law of a diffusion and a continuous martingale. The process $\{X_{t \wedge T}, t \geq 0\}$ under $\mathbb{P}^x$ is a bounded continuous martingale and, therefore, it is a time change of Brownian motion. On the interval $(0, 1)$ the diffusion coefficient for $X$ is nondegenerate and the infinitesimal generator for $X$ killed upon exiting $(0, 1)$ is

$$\mathcal{L}f(x) = \tfrac{1}{2} x^{2\alpha} f''(x)$$

with zero Dirichlet boundary conditions at endpoints 0 and 1. Therefore, $\mathbb{Q}_1^x$ is the law of a diffusion with infinitesimal generator

(2.3) $$\frac{\mathcal{L}(hf)(x)}{h(x)} = \frac{1}{2}|x|^{2\alpha}f''(x) + |x|^{2\alpha-1}f'(x).$$

Thus $\mathbb{Q}_1^x$ is the law of a time change of a three-dimensional Bessel process killed upon hitting 1. More precisely, if we let $\tau_t := \inf\{s \geq 0 : \int_0^s |X_r|^{2\alpha}\,dr \geq t\}$ for $t \geq 0$, then under $\mathbb{Q}_1^x$, the time-changed process $\{X_{\tau_t \wedge T}, t \geq 0\}$ is a three-dimensional Bessel process starting from $x$ and killed upon hitting 1.

REMARK 2.3. We will need to use the fact that if $X$ and $Y$ are two weak solutions to (1.1) and (1.3) with $X_0 \neq Y_0$ that are driven by a common Brownian motion $W$, then $X_t \neq Y_t$ for $t < \inf\{s > 0 : |X_s| + |Y_s| = 0\}$. This is well known, but we indicate the proof for the convenience of the reader. Let $X_t^x$ denote the solution to (1.2) started at $x$ and stopped at the hitting time of 0. Since $x \mapsto |x|^\alpha$ is smooth except at 0, the process $X_t^x$ is unique in the pathwise sense. Moreover we can choose versions of $X_t^x$ such that the map

is smooth on $\{x > 0, 0 \leq t < T_0^{X^x}\}$. Informally speaking, we have a flow; see [3]. Write $b(x) = |x|^\alpha$. If $D_t^x = \frac{\partial X_t^x}{\partial x}$, then $D_t^x$ solves the equation

$$D_t^x = 1 + \int_0^t b'(X_s^x) D_s^x \, dW_s, \qquad t < T_0^{X^x}.$$

This is a linear equation with the unique solution

$$D_t^x = \exp\left(\int_0^t b'(X_s^x) \, dW_s - \tfrac{1}{2}\int_0^t (b'(X_s^x))^2 \, ds\right), \qquad t < T_0^{X^x}.$$

Therefore, except for a null set, $D_t^x$ is strictly positive on $\{x > 0, 0 \leq t < T_0^{X^x}\}$. So when $X_0 > 0$ and $Y_0 < X_0$, this implies $Y_t < X_t$ for $t < T_0^X$, and similarly for other orderings of the three points $\{0, X_0, Y_0\}$.

Throughout we let the letter $c$ with or without subscripts denote constants whose exact value is unimportant and may change from line to line.

**3. The chasing phenomenon.** Recall that $X$ and $Y$ denote solutions to (1.1) and (1.3) driven by the same Brownian motion. We begin by showing that if $X$ is conditioned to hit 1 before 0, then $Y$ will "chase" after $X$ and will be close to 1 when $X$ hits 1, provided $|Y_0| \leq X_0$ and $X_0$ is small.

LEMMA 3.1. *For $x \in (0, 1)$, we have*

$$(3.1) \qquad \mathbb{E}^x\left[\int_0^{T_1^X} X_s^{2\alpha - 2} \, ds \bigg| T_1^X < T_0^X\right] = -2\log x.$$

PROOF. The Green function for a three-dimensional Bessel process starting from $x$ and killed upon hitting 1 is

$$G(x, y) = \begin{cases} 2y^2\left(\dfrac{1}{y} - 1\right), & y \in (x, 1), \\ 2y^2\left(\dfrac{1}{x} - 1\right), & y \in (0, x). \end{cases}$$

The Green function for the process $X$ under $\mathbb{Q}_1^x$ is then $G(x,y)y^{-2\alpha}$. Therefore

$$\begin{aligned}
\mathbb{E}^x\left[\int_0^{T_1} X_s^{2\alpha - 2} \, ds \bigg| T_1 < T_0\right] &= \int_0^1 y^{2\alpha - 2} G(x, y) y^{-2\alpha} \, dy \\
&= \int_0^x 2\left(\frac{1}{x} - 1\right) dy + \int_x^1 2\left(\frac{1}{y} - 1\right) dy \\
&= -2\log x. \qquad \square
\end{aligned}$$

Assume for a moment that $X_0 > 0$ and $Y_0 > 0$. Applying Itô's formula, we have for $t < T_0^X$,

$$
\begin{aligned}
dX_t^{1-\alpha} &= (1-\alpha)X_t^{-\alpha}\,dX_t - \tfrac{1}{2}\alpha(1-\alpha)X_t^{-\alpha-1}\,d\langle X\rangle_t \\
&= (1-\alpha)\,dW_t - \tfrac{1}{2}\alpha(1-\alpha)X_t^{\alpha-1}\,dt.
\end{aligned}
\tag{3.2}
$$

Similarly, for $t < T_0^Y$,

$$
dY_t^{1-\alpha} = (1-\alpha)\,dW_t - \tfrac{1}{2}\alpha(1-\alpha)Y_t^{\alpha-1}\,dt.
\tag{3.3}
$$

Let

$$R_t = X_t^{1-\alpha} - Y_t^{1-\alpha}.$$

Then (3.2) and (3.3) imply that for $t < T_0^X \wedge T_0^Y$,

$$
\begin{aligned}
dR_t &= \frac{1}{2}\alpha(1-\alpha)\left(\frac{1}{Y_t^{1-\alpha}} - \frac{1}{X_t^{1-\alpha}}\right)dt \\
&= \frac{1}{2}\alpha(1-\alpha)R_t X_t^{\alpha-1}Y_t^{\alpha-1}\,dt.
\end{aligned}
$$

This is an ordinary differential equation and we obtain for $t < T_0^X \wedge T_0^Y$

$$
R_t = R_0 \exp\left(\tfrac{1}{2}\alpha(1-\alpha)\int_0^t X_s^{\alpha-1}Y_s^{\alpha-1}\,ds\right).
\tag{3.4}
$$

LEMMA 3.2. *For every $\delta \in (0,1)$ there exists $\kappa_0 \in (0,1)$ (depending only on $\delta$) such that if*

$$1/2 > X_0 > Y_0 > (1-\kappa_0)X_0 > 0,$$

*then*

$$\mathbb{P}^{X_0,Y_0}(Y_t > (1-\delta)X_t \text{ for all } t \leq T_1^X \mid T_1^X < T_0^X) \geq 1-\delta.$$

Recall that $\mathbb{P}^{x,y}$ denotes any element of $\mathcal{P}(x,y)$, and so the above lemma asserts the estimate for every element of $\mathcal{P}(X_0, Y_0)$.

PROOF OF LEMMA 3.2. It suffices to consider the case when $X_0 \in (0, 1/2)$ is deterministic, say $X_0 = x_0$. Choose $j_0$ so that $2^{-j_0} \leq x_0 < 2^{-j_0+1}$. For notational simplicity, let

$$\sigma_k = T_{2^{-j_0+k}}^X = \inf\{t: X_t = 2^{-j_0+k}\},$$

and define

$$\xi_k = \tfrac{1}{2}\alpha(1-\alpha)\int_{\sigma_k}^{\sigma_{k+1}} X_s^{2\alpha-2}\,ds.$$

Remark 2.2 tells us that $\{X_t, t < T\}$ under $\mathbb{Q}_\infty^x$ is a time change of a three-dimensional Bessel process. By the pseudo-strong Markov property of $X$ and scaling, under $\mathbb{Q}_\infty^{x_0}$ the $\{\xi_k, k \geq 1\}$ are i.i.d. random variables, having the same distribution as $\frac{1}{2}\alpha(1-\alpha) \int_0^{T_1^X} X_s^{2\alpha-2} ds$ under $\mathbb{Q}_\infty^{1/2}$. By (3.1) we have

$$(3.5) \quad \mathbb{E}^{1/2}\left[\frac{1}{2}\alpha(1-\alpha)\int_0^{T_1^X} X_s^{2\alpha-2} ds \Big| T_1^X < T_0^X\right] = \alpha(1-\alpha)\log 2.$$

Set $\gamma = 9/8$ and define

$$(3.6) \quad B_1 = \left\{\sum_{k=0}^N \xi_k \leq c_1 + N\alpha(1-\alpha)\gamma\log 2 \text{ for every } N \in [1, j_0 - 1]\right\}.$$

We claim there exists $c_1$ such that

$$\mathbb{Q}_1^{x_0}(B_1) \geq 1 - \delta/2.$$

Indeed, by the strong law of large numbers there exists $n_0$ such that

$$\mathbb{Q}_\infty^{x_0}\left(\sum_{k=0}^N \xi_k \leq N\alpha(1-\alpha)\gamma\log 2 \text{ for all } N \geq n_0\right) \geq 1 - \delta/4.$$

We then choose $c_1$ sufficiently large to make $\mathbb{Q}_\infty^{x_0}(B_1) > 1 - \delta/2$, and note that $\mathbb{Q}_1^{x_0}(B_1) = \mathbb{Q}_\infty^{x_0}(B_1)$.

Without loss of generality, assume that $0 < \delta < 1 - (7/8)^{1/(1-\alpha)}$. Choose $\beta \in (0, 1/8)$ such that

$$(3.7) \quad (1-\beta)^{1/(1-\alpha)} \geq 1 - \delta.$$

Then on the event $B_1$ we have for $t \in [\sigma_{N-1}, \sigma_N)$

$$(3.8) \quad \begin{aligned} X_0^{1-\alpha} \exp&\left(\frac{1}{2}\frac{\alpha(1-\alpha)}{1-\beta}\int_0^t X_s^{2-2\alpha} ds\right) \\ &\leq 2^{(-j_0+1)(1-\alpha)} \exp\left(\frac{1}{2}\frac{\alpha(1-\alpha)}{1-\beta}\int_0^{\sigma_N} X_s^{2-2\alpha} ds\right) \\ &\leq 2^{(-j_0+1)(1-\alpha)} e^{c_1/(1-\beta)} e^{N\alpha(1-\alpha)\gamma\log 2/(1-\beta)} \\ &= c_2 2^{-j_0(1-\alpha)} (2^N)^{\alpha(1-\alpha)\gamma/(1-\beta)}. \end{aligned}$$

By Remark 2.2, the process $\{X_t, t < T\}$ under $\mathbb{Q}_1^x$ is a time change of a three-dimensional Bessel process, whereas a three-dimensional Bessel process has the same distribution as the modulus of a three-dimensional Brownian motion. By Proposition I.5.8(b)(iii) of [2] the probability that $X$ under $\mathbb{Q}_1^{2^{-j_0+k}}$ will ever hit $c_3 2^{-j_0+k}/2^{k/8} = c_3 2^{-j_0+7k/8}$ is $c_3 2^{-k/8}$. Summing over $k$ from 1 to $\infty$ and taking $c_3$ small enough, there is $\mathbb{Q}_1^{x_0}$ probability at most $\delta/2$ that $X_t$ gets below $c_3 2^{-j_0+7k/8}$ between times $\sigma_k$ and $\sigma_{k+1}$ for some $k \geq 1$. Let

$$(3.9) \qquad B_2 = \left\{ \inf_{\sigma_k \leq t < \sigma_{k+1}} X_t \geq c_3 2^{-j_0+7k/8} \text{ for } k = 1, 2, \ldots, j_0 - 1 \right\},$$

and $B = B_1 \cap B_2$. Since $0 < \alpha < 1/2$ and $0 < \beta < 1/8$, except for the event $B^c$ of $\mathbb{Q}_1^{x_0}$-probability at most $\delta$, we have from (3.8) and (3.9) that for $t \in [\sigma_{N-1} \sigma_N)$,

$$\begin{aligned}
X_0^{1-\alpha} \exp&\left( \frac{1}{2} \frac{\alpha(1-\alpha)}{1-\beta} \int_0^t X_s^{2-2\alpha} \, ds \right) \\
&\leq c_2 2^{-j_0(1-\alpha)} (2^N)^{\alpha(1-\alpha)\gamma/(1-\beta)} \\
(3.10) \qquad &\leq c_2 (2^{-j_0 + 9N\alpha/(8(1-\beta))})^{1-\alpha} \\
&\leq c_2 (2^{-j_0 + 7N/8})^{1-\alpha} \\
&\leq c_4 X_t^{1-\alpha}.
\end{aligned}$$

Define

$$S = \inf\{t > 0 : Y_t^{1-\alpha} \leq (1-\beta) X_t^{1-\alpha}\}.$$

Under $\mathbb{Q}_1^{x_0}$, the process $\{X_t, t \leq T\}$ never hits 0. Since $Y_0 > 0$, then the process $Y$ cannot hit 0 before time $S$. Choose $\kappa = \beta/(2c_4) \wedge \frac{1}{2}$ and let us define $\kappa_0 = 1 - (1-\kappa)^{1/(1-\alpha)}$. Then the condition $Y_0 > (1-\kappa_0)X_0 > 0$ implies that $Y_0^{1-\alpha} > \kappa X_0^{1-\alpha}$. It follows from (3.4) and (3.10) that on the event $B$, for every $t \leq S$,

$$\begin{aligned}
R_t &\leq R_0 \exp\left( \frac{\alpha(1-\alpha)}{2(1-\beta)} \int_0^t X_s^{2\alpha-2} \, ds \right) \\
&\leq \kappa X_0^{1-\alpha} \exp\left( \frac{\alpha(1-\alpha)}{2(1-\beta)} \int_0^t X_s^{2\alpha-2} \, ds \right) \\
&\leq \kappa c_4 X_t^{1-\alpha} \leq \frac{\beta}{2} X_t^{1-\alpha}.
\end{aligned}$$

Since $R_S = \beta X_S^{1-\alpha}$ on $\{S < \infty\}$, we conclude from above that $S = \infty$ on $B$ except for a $\mathbb{Q}_1^{x_0}$ null set.

We have thus shown that under the condition $1/2 > X_0 > Y_0 > (1-\kappa_0)X_0 > 0$, on the event $B$ we have

$$(3.11) \qquad Y_t^{1-\alpha} > (1-\beta) X_t^{1-\alpha} \qquad \text{for every } t \geq 0.$$

It follows from the above and (3.7) that

$$\mathbb{P}^{X_0, Y_0}(Y_t > (1-\delta) X_t \text{ for all } t \leq T_1^X | T_1^X < T_0^X) \geq 1 - \delta.$$

This proves the lemma. $\square$

REMARK 3.3. Suppose $X$ and $Y$ are weak solutions to (1.1) and (1.3) driven by the same Brownian motion $W$. Since the process $X$ under $\mathbb{Q}_1^z$ has infinitesimal generator given by (2.3), there exists a $\mathbb{Q}_1$-Brownian motion $\widetilde{W}_t$ such that

$$dX_t = X_t^\alpha \, d\widetilde{W}_t + X_t^{2\alpha-1} \, dt. \tag{3.12}$$

We also have $dX_t = X_t^\alpha \, dW_t$, so

$$dW_t = d\widetilde{W}_t + X_t^{\alpha-1} \, dt. \tag{3.13}$$

Therefore,

$$dY_t = Y_t^\alpha \, d\widetilde{W}_t + X_t^{\alpha-1} Y_t^\alpha \, dt \tag{3.14}$$

for $t < T_0^X$.

We thus see that under $\mathbb{Q}_1^{x,y}$, $X$ and $Y$ solve the SDEs (3.12) and (3.14). The discussion in Remark 2.1 then shows us that the pseudo-strong Markov property holds if $\overline{\mathbb{P}} \in \mathcal{Q}_1(x,y)$. That is, if $S$ is a stopping time, $\overline{\mathbb{P}}_S$ is defined by $\overline{\mathbb{P}}_S(A) = \overline{\mathbb{P}}(A \circ \theta_S)$ for $A \in \mathcal{F}_\infty$, and $\mathbb{O}_S$ is a regular conditional probability for $\overline{\mathbb{P}}_S$ given $\mathcal{F}_S$, then $\mathbb{O}_S \in \mathcal{Q}_1(X_S, Y_S)$ almost surely. A similar argument shows that the pseudo-strong Markov property holds for $\overline{\mathbb{P}} \in \mathcal{Q}_i(x,y)$ for $i = 0$ and $i = \infty$.

LEMMA 3.4. *For every $\delta \in (0,1)$ there exist $\rho > 0$ and $K \geq 1$ such that for $X_0 = x \in (0, 2^{-K})$, $Y_0 = y \in [-x, x)$ and $S = \inf\{t > 0 : X_t = x2^K\}$,*

$$\mathbb{P}^{x,y}(Y_S > (1-\delta)X_S | T_1^X < T_0^X) \geq \rho.$$

PROOF. Putting (3.13) into (3.2) and (3.3), we have for $t < T_0^X$

$$\begin{aligned}
dX_t^{1-\alpha} &= (1-\alpha) \, d\widetilde{W}_t + (1-\alpha)X_t^{\alpha-1} \, dt - \frac{1}{2}\alpha(1-\alpha)X_t^{\alpha-1} \, dt \\
&= (1-\alpha) \, d\widetilde{W}_t + (1-\alpha)\left(1 - \frac{\alpha}{2}\right) X_t^{\alpha-1} \, dt.
\end{aligned} \tag{3.15}$$

For $Y_0 \neq 0$, applying Itô's formula to $Y_t$ and using (3.14), we have for $t < T_0^Y$,

$$\begin{aligned}
d|Y_t|^{1-\alpha} &= (1-\alpha)\operatorname{sgn}(Y_t) \, d\widetilde{W}_t + (1-\alpha)\operatorname{sgn}(Y_t)X_t^{\alpha-1} \, dt \\
&\quad - \tfrac{1}{2}\alpha(1-\alpha)|Y_t|^{\alpha-1} \, dt.
\end{aligned} \tag{3.16}$$

We first consider the case where $Y_0 < 0$. Let

$$A_1 = \left\{ \sup_{0 \leq s \leq x^{2-2\alpha}} \widetilde{W}_s > 3x^{1-\alpha} \text{ and } \inf_{0 \leq s \leq x^{2-2\alpha}} \widetilde{W}_s > -x^{1-\alpha}/2 \right\}.$$

By Brownian scaling, $\mathbb{Q}_1^z(A_1) > c_2$ with $c_2$ independent of $x$ and $z$. It follows from (3.16) that for $t < T_0^Y$,

$$|Y_t|^{1-\alpha} \leq |Y_0|^{1-\alpha} - (1-\alpha)\widetilde{W}_t$$
$$\leq x^{1-\alpha} - (1-\alpha)\widetilde{W}_t.$$

If $A_1$ holds, then

$$\inf_{0 \leq t \leq x^{2-2\alpha}} (x^{1-\alpha} - (1-\alpha)\widetilde{W}_t) \leq x^{1-\alpha} - 3(1-\alpha)x^{1-\alpha} \leq -\tfrac{1}{2}x^{1-\alpha} < 0.$$

This implies that $T_0^Y < x^{2-2\alpha}$ on $A_1$. On the other hand, on the event $A_1$,

$$\inf_{0 \leq t \leq x^{2-2\alpha}} (x^{1-\alpha} + (1-\alpha)\widetilde{W}_t) \geq x^{1-\alpha} - \frac{1-\alpha}{2}x^{1-\alpha} \geq \frac{1}{2}x^{1-\alpha} > 0.$$

So by (3.15) on $A_1$,

$$T_0^X > x^{2-2\alpha} > T_0^Y \quad \text{and} \quad \inf_{0 \leq t \leq x^{2-2\alpha}} X_t^{1-\alpha} \geq \tfrac{1}{2}x^{1-\alpha}.$$

Note also that on the event $A_1$, by (3.15) again,

$$\sup_{0 \leq s \leq x^{2-2\alpha}} X_s^{1-\alpha} \leq x^{1-\alpha} + (1-\alpha)3x^{1-\alpha} + (1-\alpha)\left(1 - \frac{\alpha}{2}\right)\int_0^{x^{2-2\alpha}} 2x^{\alpha-1}\,ds$$
$$\leq 6x^{1-\alpha}.$$

So in particular we have that

$$2^{-1/(1-\alpha)}x \leq X_{T_0^Y} \leq 6^{1/(1-\alpha)}x \qquad \text{on } A_1.$$

Using the pseudo-strong Markov property at the time $\inf\{t \geq 0 : Y_t \geq 0\}$, it thus suffices to prove that the conclusion of the lemma holds under the following assumptions: $0 \leq Y_0 < X_0$ and

$$2^{-1/(1-\alpha)}x \leq X_0 \leq 6^{1/(1-\alpha)}x.$$

As this is our first nontrivial use of the pseudo-strong Markov property, we explain in detail how it is used. Recall that $T := T_0^X \wedge T_1^X$. We use the pseudo-strong Markov property at the time $V = \inf\{t \geq 0 : Y_t \geq 0\} \wedge T$. Let $\mathbb{O}_V$ be a regular conditional probability for the law of $\{(X_{V+t}, Y_{V+t}); 0 \leq t \leq T \circ \theta_V\}$ under $\mathbb{Q}_1^{x,y}$ conditional on $\mathcal{F}_V$. This means that for each $A \in \mathcal{F}_\infty$ the random variable $\mathbb{O}_V(\cdot, A)$ is $\mathcal{F}_V$-measurable, $\mathbb{O}_V(\omega, \cdot)$ is a probability on $\mathcal{F}_T$ for each $\omega$, and for each $B \in \mathcal{F}_V$ and each $A$ of the form

$$A = \{(X_{t_1 \wedge T}, Y_{t_1 \wedge T}) \in C_1, \ldots, (X_{t_n \wedge T}, Y_{t_n \wedge T}) \in C_n\}$$

for positive reals $t_1 \leq \cdots \leq t_n$ and $C_1, \ldots, C_n$ Borel subsets of $\mathbb{R}^2$,

$$\mathbb{E}_{\mathbb{Q}_1^{x,y}}[\mathbb{O}_V(\cdot, A); B]$$
$$= \mathbb{Q}_1^{x,y}(\{(X_{V+t_1}, Y_{V+t_1}) \in C_1, \ldots, (X_{V+t_n}, Y_{V+t_n}) \in C_n\} \cap B).$$

Here $\mathbb{E}_{\mathbb{Q}_1^{x,y}}$ is the expectation with respect to $\mathbb{Q}_1^{x,y}$. As in the discussions in Remarks 2.1 and 3.3, for $\mathbb{Q}_1^{x,y}$-almost every $\omega \in \{V < T\}$, the law of $\{(X_t, Y_t); t \leq T\}$ under $\mathbb{O}_V(\omega, \cdot)$ is that of the law of a weak solution $(\widetilde{X}, \widetilde{Y})$ to (3.12) and (3.14) with $(\widetilde{X}_0(\omega'), \widetilde{Y}_0(\omega')) = (X_V(\omega), Y_V(\omega))$ for $\mathbb{O}_V(\omega, \cdot)$-almost every $\omega'$. If we let $S' = \inf\{t \geq V : X_t = x2^K\}$, then

$$\mathbb{Q}_1^{x,y}(Y_S > (1-\delta)X_S) \geq \mathbb{Q}_1^{x,y}(\{Y_{S'} > (1-\delta)X_{S'}\} \cap A_1)$$
$$= \mathbb{E}_{\mathbb{Q}_1^{x,y}}[\mathbb{O}_V(\cdot, Y_S > (1-\delta)X_S); A_1].$$

Since we have shown $\mathbb{Q}_1^{x,y}(A_1)$ is bounded below, to prove the lemma it suffices to find a lower bound on $\mathbb{O}_V(\cdot, Y_S > (1-\delta)X_S)$. By the remarks above and the fact that $X_V \in [2^{1/(1-\alpha)}x, 6^{1/(1-\alpha)}x]$ on the event $A_1$, this means we have to show that $\mathbb{Q}_1^{X_0, Y_0}(Y_S > (1-\delta)X_S)$ is bounded below when we have $Y_0 = 0$ and $X_0 \in [2^{1/(1-\alpha)}x, 6^{1/(1-\alpha)}x]$.

As $Y$ solves $dY_t = |Y_t|^\alpha \, dW_t$ and spends zero time at 0, the Green function for the subprocess of $Y$ starting from $y \in (-\eta, \eta)$ killed upon exiting $(-\eta, \eta)$ is

$$G(y, z) = \begin{cases} |z|^{-2\alpha}(y + \eta)(\eta - z)/\eta, & -\eta < y < z < \eta, \\ |z|^{-2\alpha}(z + \eta)(\eta - y)/\eta, & -\eta < z < y < \eta. \end{cases}$$

By our notation convention, $T_\eta^{|Y|} = \inf\{t : |Y_t| = \eta\}$. Then for $0 \leq y < \eta$,

$$\mathbb{E}^y[T_\eta^{|Y|}] = \int_{-\eta}^{\eta} G(y, z) \, dz \leq c_3 \eta^{2-2\alpha}.$$

Let $c_4 = (1/2)^{1/(1-\alpha)}$ and $c_5 = 6^{1/(1-\alpha)}$, and so $c_4 x \leq X_0 \leq c_5 x$. By scaling, we can choose $c_6$ such that

$$\mathbb{P}^x\left(\sup_{0 \leq s \leq c_6 x^{2-2\alpha}} |X_s - x| > c_4 x/2\right) \leq \tfrac{1}{8}.$$

Choose $\eta$ so that $c_3 \eta^{2-2\alpha} = c_6 x^{2-2\alpha}/8$. Then for $0 \leq y < \eta$,

$$\mathbb{P}^y(T_\eta^{|Y|} \geq c_6 x^{2-2\alpha}) \leq \frac{\mathbb{E} T_\eta^{|Y|}}{c_6 x^{2-2\alpha}} \leq \frac{c_3 \eta^{2-2\alpha}}{c_6 x^{2-2\alpha}} = \frac{1}{8}$$

and so $\mathbb{P}^y(T_\eta^{|Y|} < c_6 x^{2-2\alpha}) \geq 7/8$. By symmetry, if $y \in [0, \eta)$ then $Y$ starting at $y$ will exit $(-\eta, \eta)$ through $\eta$ with probability at least $1/2$. Hence

$$\mathbb{P}^y(T_\eta^Y < c_6 x^{2-2\alpha} | Y_0 \in [0, \eta)) \geq \tfrac{7}{8} - \tfrac{1}{2} = \tfrac{3}{8}.$$

Let

$$A_2 = \left\{\sup_{0 \leq s \leq c_6 x^{2-2\alpha}} |X_s - x| \leq c_4 x/2 \text{ and } T_\eta^Y < c_6 x^{2-2\alpha}\right\}.$$

Note that on $A_2$,

$$c_4 x/2 \leq X_t \leq (c_5 + (c_4/2))x \qquad \text{for } t \leq c_6 x^{2-2\alpha}.$$

Write $t_0 = c_6 x^{2-2\alpha}$. If $y \in [0, \eta)$, then

$$\mathbb{P}^{X_0,y}(A_2 \cap \{T_1^X < T_0^X\})$$
$$= \mathbb{P}^{X_0,y}\Big(T_\eta^Y < t_0, \sup_{0 \leq s \leq t_0} |X_s - x| \leq c_4 x/2, t_0 < T_1^X < T_0^X\Big).$$

Using the pseudo-strong Markov property of $X$ at time $t_0$, this is bounded below by

$$\mathbb{P}^{X_0,y}(A_2)\mathbb{P}^{c_4 x/2}(T_1^X < T_0^X) \geq (\tfrac{3}{8} - \tfrac{1}{8})(c_4 x/2) = c_4 x/8.$$

Since $\mathbb{P}^{X_0}(T_1^X < T_0^X) \leq c_5 x$, we conclude that

$$\mathbb{Q}_1^{X_0,y}(A_2) \geq \frac{c_4}{8c_5}.$$

Applying the pseudo-strong Markov property at the stopping time $\inf\{t \geq 0 : Y_t \geq \eta\}$, we may thus assume that

$$c_7 x \leq X_0 \leq c_8 x \quad \text{and} \quad Y_0 \geq c_9 x,$$

where $c_7 = c_4/2 = 2^{-1-1/(1-\alpha)}$, $c_8 = c_5 + (c_4/2)$ and $c_9 = \tfrac{1}{2} \wedge (\tfrac{c_6}{8c_3})^{1/(2-2\alpha)}$.

In view of Lemma 3.2 and the pseudo-strong Markov property, it suffices to show that there is positive $\mathbb{Q}_1^{X_0,Y_0}$-probability that there exists $t < S$ such that

$$(3.17) \qquad Y_t \geq (1 - \kappa_0) X_t > 0$$

where $\kappa_0$ is the constant in Lemma 3.2 corresponding to $\delta$. Let

$$a = \frac{c_9^{1-\alpha}}{8(1-\alpha)}, \qquad c_{10} = c_7^{1-\alpha} - \frac{1}{8} c_9^{1-\alpha}, \qquad \gamma = \Big(\frac{c_9^{2-2\alpha}}{8\alpha(1-\alpha)}\Big) \wedge (c_9^{1-\alpha} c_{10} \log 2)$$

and

$$b = \frac{4 c_8^{1-\alpha}}{1 - (1-\kappa_0)^{1-\alpha}}.$$

Define

$$A_3 = \Big\{\sup_{0 \leq s \leq \gamma x^{2-2\alpha}} \widetilde{W}_s \in (b x^{1-\alpha}, 2 b x^{1-\alpha}) \text{ and } \inf_{0 \leq s \leq \gamma x^{2-2\alpha}} \widetilde{W}_s > -a x^{1-\alpha}\Big\}.$$

By Brownian scaling, $\mathbb{Q}_1^z(A_3)$ is bounded below by a positive constant independent of $x$ and $z$. Let $U = \inf\{t : Y_t \leq c_9 2^{-1/(1-\alpha)} x\}$. By (3.16), on the event $A_3$,

$$Y^{1-\alpha}_{U \wedge \gamma x^{2-2\alpha}} \geq Y_0^{1-\alpha} + (1-\alpha)\widetilde{W}_{U \wedge \gamma x^{2-2\alpha}} - \frac{1}{2}\alpha(1-\alpha)\int_0^{U \wedge \gamma x^{2-2\alpha}} Y_s^{\alpha-1}\,ds$$

$$\geq (c_9 x)^{1-\alpha} - (1-\alpha)ax^{1-\alpha} - \frac{\alpha(1-\alpha)\gamma}{c_9^{1-\alpha}}x^{1-\alpha}$$

$$\geq 3c_9^{1-\alpha}x^{1-\alpha}/4.$$

We conclude that on the event $A_3$ we have $U \geq \gamma x^{2-2\alpha}$. We have by (3.15) that on the event $A_3$,

$$\inf_{0 \leq t \leq \gamma x^{2-2\alpha}} X_t^{1-\alpha} \geq \inf_{0 \leq t \leq \gamma x^{2-2\alpha}} (X_0^{1-\alpha} - (1-\alpha)\widetilde{W}_t)$$

$$\geq (c_7 x)^{1-\alpha} - (1-\alpha)ax^{1-\alpha}$$

$$= c_{10}x^{1-\alpha}.$$

Using (3.4) we then have on $A_3$ for $t \leq \gamma x^{2-2\alpha}$,

$$R_t \leq R_0 \exp\left(\frac{1}{2}\alpha(1-\alpha)\int_0^{\gamma x^{2-2\alpha}} 2c_{10}^{-1}x^{\alpha-1}(c_9 x)^{\alpha-1}\,dt\right)$$

$$\leq X_0^{1-\alpha}\exp(\alpha(1-\alpha)\gamma c_{10}^{-1}c_9^{\alpha-1})$$

$$\leq 2c_8^{1-\alpha}x^{1-\alpha}.$$

On the other hand, on $A_3$ there exists some $t \leq \gamma x^{2-2\alpha}$ such that

(3.18) $$X_t^{1-\alpha} \geq (1-\alpha)\widetilde{W}_t \geq bx^{1-\alpha}/2.$$

Therefore on $A_3$ there exists some $t \leq \gamma x^{2-2\alpha}$ such that

$$R_t \leq 2c_8^{1-\alpha}x^{1-\alpha} \leq (1-(1-\kappa_0)^{1-\alpha})X_t^{1-\alpha}.$$

Thus

(3.19) $$\inf\{t \geq 0 : Y_t \geq (1-\kappa_0)X_t\} \leq \gamma x^{2-2\alpha} \quad \text{on } A_3.$$

It follows from (3.15) and (3.18) that on $A_3$,

$$\sup_{0 \leq t \leq \gamma x^{2-2\alpha}} X_t^{1-\alpha}$$

$$\leq \sup_{0 \leq t \leq \gamma x^{2-2\alpha}} \left(X_0^{1-\alpha} + (1-\alpha)\widetilde{W}_t + \int_0^t (1-\alpha)\left(1-\frac{\alpha}{2}\right)|X_s|^{\alpha-1}\,ds\right)$$

$$\leq (c_8 x)^{1-\alpha} + 2(1-\alpha)bx^{1-\alpha} + (1-\alpha)\left(1-\frac{\alpha}{2}\right)\gamma x^{2-2\alpha}\frac{2}{b}x^{1-\alpha}$$

$$\leq c_{11}^{1-\alpha}x^{1-\alpha}.$$

That is,
$$\sup_{0 \leq t \leq \gamma x^{2-2\alpha}} X_t \leq c_{11} x \qquad \text{on } A_3.$$

Now take $K \geq 1$ large enough so that $2^K > c_{11}$. Then for every $x \in (0, 2^{-K})$,
$$S := \inf\{t > 0 : X_t = x 2^K\} \in (\gamma x^{2-2\alpha}, T_1^X) \qquad \text{on } A_3.$$

This, together with (3.19), proves that (3.17) holds on $A_3$, which completes the proof of the lemma. $\square$

LEMMA 3.5. *Suppose that $X_0 = x \in (0, 1/4)$ and $Y_0 = y \in [-x, x)$. Let $T_{0+}^Y = \inf\{t \geq 0 : Y_t \geq 0\}$. Then for every $\lambda \geq 4x$, $T_{0+}^Y < T_\lambda^X$, $\mathbb{Q}_1^{x,y}$-a.s., and*
$$\max\{X_t, |Y_t|\} \leq 2^{1/(1-\alpha)} x < 4x \qquad \text{for } t \leq T_{0+}^Y,$$
$$-(\lambda^{1-\alpha} - X_t^{1-\alpha})^{1/(1-\alpha)} \leq Y_t < X_t \qquad \text{for } t \in (T_{0+}^Y, T_\lambda^X].$$

*In particular,*
$$|Y_t| < \lambda \qquad \text{for every } t \leq T_\lambda^X.$$

PROOF. As $|x|^\alpha$ is Lipschitz on $(0, \infty)$, $X_t > 0$ for every $t \geq 0$ under $\mathbb{Q}_1^{x,y}$, and $Y_0 < X_0$, it follows by a standard comparison theorem (see, e.g., Theorem I.6.2 in [3]) that

$$(3.20) \qquad Y_t \leq X_t \qquad \text{for } t < T_0^X.$$

If $Y_0 < 0$, it follows from (3.15) and (3.16) that under $\mathbb{Q}_1^{x,y}$,

$$(3.21) \qquad d(X_t^{1-\alpha} + |Y_t|^{1-\alpha}) = -\tfrac{1}{2}\alpha(1-\alpha)(X^{\alpha-1} + |Y_t|^{\alpha-1})\,dt$$

for $t < T_0^Y$. So for $t \leq T_0^Y$,
$$X_t^{1-\alpha} + |Y_t|^{1-\alpha} < |X_0|^{1-\alpha} + |Y_0|^{1-\alpha} \leq 2x^{1-\alpha},$$
and, therefore, for $t \leq T_0^Y$,
$$\max\{X_t, |Y_t|\} < 2^{1/(1-\alpha)} x < 4x.$$

Since $\lambda \geq 4x$, this implies that $T_{0+}^Y < T_\lambda^X$ and $|Y_t| < 4x \leq \lambda$ for $t < T_{0+}^Y$.

Define
$$\alpha_1 = \inf\{t > T_{0+}^Y : Y_t = -(\tfrac{1}{2}(\lambda^{1-\alpha} - X_t^{1-\alpha}))^{1/(1-\alpha)}\} \wedge T_\lambda^X,$$
$$\beta_1 = \inf\{t > \alpha_1 : Y_t = 0\} \wedge T_\lambda^X,$$

and for $n \geq 2$,
$$\alpha_n = \inf\{t > \beta_{n-1} : Y_t = -(\tfrac{1}{2}(\lambda^{1-\alpha} - X_t^{1-\alpha}))^{1/(1-\alpha)}\} \wedge T_\lambda^X,$$
$$\beta_n = \inf\{t > \alpha_n : Y_t = 0\} \wedge T_\lambda^X.$$

On $[T_{0+}^Y, \alpha_1)$,
$$-(\tfrac{1}{2}(\lambda^{1-\alpha} - X_t^{1-\alpha}))^{1/(1-\alpha)} \leq Y_t < X_t < 1.$$

On $\{\alpha_n < T_\lambda^X\}$, we have by (3.21) that $t \mapsto X_t^{1-\alpha} + |Y_t|^{1-\alpha}$ is decreasing on $[\alpha_n, \beta_n]$ and so for $t \in [\alpha_n, \beta_n]$,

$$(3.22) \quad X_t^{1-\alpha} + |Y_t|^{1-\alpha} < X_{\alpha_n}^{1-\alpha} + |Y_{\alpha_n}|^{1-\alpha} = \tfrac{1}{2}(X_{\alpha_n}^{1-\alpha} + \lambda^{1-\alpha}) < \lambda^{1-\alpha}.$$

This in particular implies that

$$(3.23) \quad \text{if } \alpha_n < T_\lambda^X, \quad \text{then } \beta_n < T_\lambda^X \quad \text{and} \quad |Y_t| < \lambda \quad \text{for } t \in [\alpha_n, \beta_n].$$

If $\alpha_n < T_\lambda^X$, then on $[\beta_n, \alpha_{n+1})$,

$$(3.24) \quad -(\tfrac{1}{2}(\lambda^{1-\alpha} - X_t^{1-\alpha}))^{1/(1-\alpha)} \leq Y_t < X_t \leq \lambda,$$

and so $|Y_t| < \lambda$ for $t \in [\beta_n, \alpha_{n+1})$. Since on $\{\alpha_n < T_\lambda^X\}$,

$$Y_{\alpha_n} = -(\tfrac{1}{2}(\lambda^{1-\alpha} - X_{\alpha_n}^{1-\alpha}))^{1/(1-\alpha)} \quad \text{and} \quad Y_{\beta_n} = 0,$$

and $Y$ has a finite number of oscillations greater than any fixed $\varepsilon > 0$ on any finite time interval, it follows that either $\alpha_n = T_\lambda^X$ for some finite $n$, or $\beta_n = T_\lambda^X$ for some finite $n$, or $\alpha_n \leq \beta_n < T_\lambda^X$ for all $n < \infty$ and $\lim_{n \to \infty} \alpha_n = \lim_{n \to \infty} \beta_n = T_\lambda^X$. If $\alpha_n = T_\lambda^X$ or $\beta_n = T_\lambda^X$ for some finite $n$ then (3.24) holds for all $t \in (T_{0+}^Y, T_\lambda^X]$ by (3.22)–(3.24). If $\alpha_n \leq \beta_n < T_\lambda^X$ for all $n < \infty$ and $\lim_{n \to \infty} \alpha_n = \lim_{n \to \infty} \beta_n = T_\lambda^X$ then (3.24) holds for all $t \in (T_{0+}^Y, T_\lambda^X)$ by (3.22)–(3.24). In this case, $Y_{T_\lambda^X} = 0$ by continuity, so (3.24) holds also for $t = T_\lambda^X$. This completes the proof of the lemma. □

LEMMA 3.6. *For every $\delta \in (0,1)$ there exists $\varepsilon > 0$ such that if $X_0 = x_0 \in (0, \varepsilon)$ and $Y_0 = y_0 \in [-x_0, x_0]$, then*

$$\mathbb{P}^{x_0, y_0}(Y_T > 1 - \delta | X_T = 1) \geq 1 - \delta,$$

*where $T = T_0^X \wedge T_1^X$.*

PROOF. When $y_0 = x_0$, then by Remark 2.3, $Y_t = X_t$ for $t \leq T$ and so the conclusion of the lemma holds. Thus without loss of generality, we assume now that $y_0 \in [-x_0, x_0)$. Let $\kappa_0$ be the constant from Lemma 3.2 that corresponds to $\delta/3$. Let $S_K = \inf\{t : X_t = 2^K x\}$. By Lemma 3.4 there exist $\rho > 0$ and integer $K \geq 1$ such that for any $x \in (0, 2^{-K})$,

$$(3.25) \quad \mathbb{Q}_1^{x, y_0}(Y_{S_K} > (1 - \kappa_0) X_{S_K}) \geq \rho.$$

Let integer $j_0 \geq 1$ be sufficiently large so that

$$1 - (1 - \rho)^{j_0} > \frac{1 - \delta}{1 - \delta/3}$$

and define $\varepsilon = 2^{-j_0 K - 2}$. Let
$$A = \{Y_{S_{jK}} > (1 - \kappa_0) X_{S_{jK}} \text{ for some } j = 1, 2, \ldots, j_0\}.$$
For every $0 < x < \varepsilon$, by Lemma 3.5,
$$(3.26) \quad -2^{j_0 K} x \leq Y_{S_{jK}} < X_{S_{jK}} = 2^{j_0 K} x \qquad \text{for every } j \in \{1, 2, \ldots, j_0\}.$$
Using the pseudo-strong Markov property at the stopping times $S_{jK}$'s recursively, we have from (3.25) that $\mathbb{Q}_1^{x_0, y_0}(A^c) \leq (1 - \rho)^{j_0}$ and so
$$\mathbb{Q}_1^{x_0, y_0}(A) \geq 1 - (1 - \rho)^{j_0}.$$
On the other hand, $A$ can be written as the disjoint union of $A_j$, $j = 1, \ldots, j_0$, where
$$A_j = \{Y_{S_{kK}} \leq (1 - \kappa_0) X_{S_{kK}} \text{ for } k \leq j - 1 \text{ and } Y_{S_{jK}} > (1 - \kappa_0) X_{S_{jK}}\}.$$
Because of (3.26), we can use the pseudo-strong Markov property at the stopping time $S_{jK}$ and apply Lemma 3.2 to conclude
$$\mathbb{Q}_1^{x_0, y_0}\left(\left\{Y_T > 1 - \frac{\delta}{3}\right\} \cap A_j\right) \geq \left(1 - \frac{\delta}{3}\right) \mathbb{Q}_1^{x_0, y_0}(A_j).$$
Therefore
$$\begin{aligned}
\mathbb{Q}_1^{x_0, y_0}\left(Y_T > 1 - \frac{\delta}{3}\right) &\geq \mathbb{Q}_1^{x_0, y_0}\left(\left\{Y_T > 1 - \frac{\delta}{3}\right\} \cap A\right) \\
&= \sum_{j=1}^{j_0} \mathbb{Q}_1^{x_0, y_0}\left(\left\{Y_T > 1 - \frac{\delta}{3}\right\} \cap A_j\right) \\
&\geq \sum_{j=1}^{j_0} \left(1 - \frac{\delta}{3}\right) \mathbb{Q}_1^{x_0, y_0}(A_j) \\
&= \left(1 - \frac{\delta}{3}\right) \mathbb{Q}_1^{x_0, y_0}(A) \geq 1 - \delta. \qquad \square
\end{aligned}$$

COROLLARY 3.7. *For every $\delta \in (0, 1)$ there exists $r_0 > 0$ such that if $X_0 \in (0, r_0)$ and $Y_0 \in [-X_0, X_0]$, then*
$$\mathbb{E}^{X_0, Y_0}[Y_T | X_T = 1] \geq 1 - \delta.$$

PROOF. When $Y_0 = X_0$, we know from Remark 2.3 that $Y_t = X_t$ for $t \leq T$ and the conclusion of the corollary holds trivially. Thus we now assume that $Y_0 \in [-X_0, X_0)$. Let $\varepsilon$ be the constant in Lemma 3.6 that corresponds to $\delta/3$. Let $r_0 = \varepsilon \wedge \frac{1}{4}$. Then by Lemma 3.5, $|Y_T| \leq 1$, $\mathbb{P}^{X_0, Y_0}$-a.s. It now follows from Lemma 3.6 that
$$\mathbb{E}^{X_0, Y_0}[Y_T | T_1^X < T_0^X] \geq \left(1 - \frac{\delta}{3}\right)^2 - \frac{\delta}{3} > 1 - \delta. \qquad \square$$

LEMMA 3.8. *Let $X$ and $Y$ be weak solutions to (1.1) and (1.3) driven by the same Brownian motion $W$ with $X_0 = x \in (0, 1/4)$ and $Y_0 = y \in [-x, x)$. Then for $\lambda \geq 4x$ we have*

$$\mathbb{P}^{x,y}(|Y_t| \leq \lambda \text{ for } t \leq T_\lambda^X \wedge T_0^X | T_0^X < T_1^X) = 1.$$

*It follows that*

$$\mathbb{P}^{x,y}(|Y_t| < 1 \text{ for } t \leq T_0^X | T_0^X < T_1^X) = 1.$$

PROOF. Consider any $a \in (0, x)$ and $b \in (x, 1)$. By Girsanov's theorem, the distributions $\mathbb{Q}_0^{x,y}$ and $\mathbb{Q}_1^{x,y}$ are mutually absolutely continuous on $[0, T_a^X \wedge T_b^X]$. Hence, by Lemma 3.5,

$$\mathbb{P}^{x,y}(|Y_t| < \lambda \text{ for } t \leq T_\lambda^X \wedge T_0^X \wedge T_a^X \wedge T_b^X | T_0^X < T_1^X) = 1.$$

Letting $a \downarrow 0$ and $b \uparrow 1$, we obtain the first result. The second result holds as on $\{T_0^X < T_1^X\}$, there is some (random) $\lambda \in (0, 1)$ such that $T_0^X = T_0^X \wedge T_\lambda^X$. □

COROLLARY 3.9. *Let $X$ and $Y$ be weak solutions to (1.1) and (1.3) driven by the same Brownian motion $W$ with $X_0 = x \in (0, 1/4)$ and $Y_0 = y \in [-x, x)$. Then for $\lambda \geq 4x$ we have*

$$\mathbb{P}^{x,y}(|Y_t| \leq \lambda \text{ for } t \leq T_0^X \wedge T_\lambda^X) = 1.$$

*In particular,*

$$\mathbb{P}^{x,y}(|Y_t| < 1 \text{ for } t \leq T_0^X \wedge T_1^X) = 1.$$

PROOF. The corollary follows immediately from Lemmas 3.5 and 3.8. □

LEMMA 3.10. *For every $\delta \in (0, 1)$, there is $r_0 > 0$ such that if $X_0 = x \in (0, r_0)$ and $Y_0 = y \in [-x, 0]$, we have $-1 \leq Y_T \leq 0$, $\mathbb{P}^{x,y}$-a.s., on $\{X_T = 0\}$, and*

(3.27) $$\mathbb{E}^{x,y}[-Y_T | X_T = 0] \geq (1 - \delta/2)x.$$

*Moreover, for every $\delta \in (0, 1)$, there exists a function $\psi : (0, 1) \to \mathbb{R}$ such that $\lim_{x \downarrow 0} \psi(x) = 0$, and for $x \in (0, r_0)$ and $y \in [-x, 0]$,*

(3.28) $$\mathbb{E}^{x,y}[-Y_T \mathbf{1}_{\{|Y_T| \leq \psi(x)\}} | X_T = 0] \geq (1 - \delta)x,$$

*or equivalently,*

(3.29) $$\mathbb{E}^{x,y}[-Y_T \mathbf{1}_{\{|Y_T| \leq \psi(x)\}} \mathbf{1}_{\{X_T = 0\}}] \geq (1 - \delta)x(1 - x).$$

PROOF. Let $\delta \in (0,1)$. By Corollary 3.7, there is $r_0 \in (0,1/4)$ so that

$$\mathbb{E}^{x,y}[Y_T|X_T=1] \geq 1 - \frac{\delta}{2}$$

whenever $X_0 = x \in (0, r_0]$ and $Y_0 = y \in [-x, x]$. Since by Corollary 3.9, $\sup_{s \leq T} |Y_s| \leq 1$ under $\mathbb{P}^{x,y}$ we have $\mathbb{E}^{x,y} Y_T = \mathbb{E}^{x,y} Y_0$ by optional stopping. The process $X$ is a continuous local martingale so $\mathbb{P}^x(X_T = 1) = x$ and $\mathbb{P}^x(X_T = 0) = 1 - x$. Hence if $X_0 = x \in (0, r_0)$ and $Y_0 = y \in [-x, 0]$,

$$\begin{aligned}
0 &\geq \mathbb{E}^{x,y} Y_0 = \mathbb{E}^{x,y}[Y_T] \\
&= \mathbb{E}^{x,y}[Y_T|X_T=1]\mathbb{P}(X_T=1) + \mathbb{E}^{x,y}[Y_T|X_T=0]\mathbb{P}(X_T=0) \\
&\geq \left(1 - \frac{\delta}{2}\right) x + \mathbb{E}^{x,y}[Y_T|X_T=0](1-x).
\end{aligned}$$

It follows that

$$(3.30) \qquad \mathbb{E}^{x,y}[-Y_T|X_T=0] \geq \frac{(1-\delta/2)x}{1-x} \geq \left(1 - \frac{\delta}{2}\right) x.$$

That $-1 \leq Y_T \leq 0$, $\mathbb{P}^{x,y}$-a.s. on $\{T_0^X < T_1^X\}$, is a consequence of Lemma 3.8 and Remark 2.3.

Consider $y_0 \in (0, 1/2)$. Suppose that $x \in (0, y_0/4)$ and note that, by Lemma 3.8 and scaling, for $y \in [-x, 0]$,

$$(3.31) \quad \mathbb{P}^{x,y}(Y_t \in [-y_0, y_0] \text{ for every } t \leq T_0^X \wedge T_{y_0}^X | T_0^X < T_{y_0}^X) = 1.$$

Recall that $\mathbb{Q}_0^x$ satisfies $\mathbb{Q}_0^x(A) = \mathbb{P}^x(A|X_T=0)$ and $\mathbb{Q}_0^x$ is derived from $\mathbb{P}$ by Doob's $h$-transform with $h(x) = 1 - x$. As in the proof of (2.3) we can show that the process $X$ under $\mathbb{Q}_0^{\cdot}$ has generator

$$(3.32) \qquad \widehat{\mathcal{L}} f = \frac{\mathcal{L}(hf)}{h} = \frac{|x|^{2\alpha}}{2} f'' - \frac{|x|^{2\alpha}}{1-x} f'.$$

Since $\mathbb{P}^x(T_{y_0}^X < T_0^X) = x/y_0$ for $x \in (0, y_0)$, we have for $0 < x < y_0 < 1$,

$$(3.33) \quad \mathbb{Q}^x(T_{y_0}^X < T_0^X) = \frac{h(y_0)}{h(x)} \mathbb{P}^x(T_{y_0}^X < T_0^X) = \frac{1-y_0}{1-x} \frac{x}{y_0} \leq \frac{x}{y_0}.$$

Suppose that $X$ has distribution $\mathbb{P}^x$. Then $X_T$ converges to 0 in distribution as $x \downarrow 0$. This and the fact that $Y$ is continuous imply that for some $y_1 \in (0, \delta y_0/16)$,

$$\mathbb{P}^{x,0}(Y_t \in [-y_0, y_0] \text{ for } t \leq T_0^X) \geq 1 - \frac{\delta y_0}{16} \qquad \text{for every } x \in (0, y_1].$$

For $x \in (0, y_1]$,

$$\mathbb{Q}_0^{x,0}(Y_t \in [-y_0, y_0] \text{ for } t \leq T_0^X)$$

$$\begin{aligned}
&= \mathbb{P}^{x,0}(Y_t \in [-y_0, y_0] \text{ for } t \leq T_0^X | X_T = 0) \\
&\geq \mathbb{P}^{x,0}(\{Y_t \in [-y_0, y_0] \text{ for } t \leq T_0^X\} \cap \{X_T = 0\}) \\
&\geq \mathbb{P}^{x,0}(Y_t \in [-y_0, y_0] \text{ for } t \leq T_0^X) - \mathbb{P}^{x,0}(X_T = 1) \\
&\geq 1 - \frac{\delta y_0}{16} - \frac{\delta y_0}{16} \\
&= 1 - \frac{\delta y_0}{8}.
\end{aligned} \quad (3.34)$$

Let $\widetilde{\mathbb{Q}}_0^x$ denote the distribution of $\{X_{t \wedge T_0^X}, t \geq 0\}$ when $\{X_t, t \geq 0\}$ has distribution $\mathbb{Q}_0^x$. Since the coefficients of the generator (3.32) are smooth except at 0, there exists a stochastic flow $\{X_t^y, t \geq 0\}$ driven by the same Brownian motion, such that $X_0^y = y$ for $y \in (0, 1)$, and the distribution of $\{X_t^y, t \geq 0\}$ is $\widetilde{\mathbb{Q}}_0^y$. So with probability 1, for every $t \geq 0$, $v < z$ implies that $X_t^v \leq X_t^z$ and the function $y \to X_t^y$ is continuous.

Let $\widetilde{\mathbb{Q}}_0^{x,y}$ denote the distribution of $\{(X_{t \wedge T_0^X}, Y_{t \wedge T_0^Y}), t \geq 0\}$ when $\{(X_t, Y_t), t \geq 0\}$ has distribution $\mathbb{Q}_0^{x,y}$. The above remarks about the stochastic flow imply that if $\{(X_t, Y_t), t \geq 0\}$ has distribution $\widetilde{\mathbb{Q}}_0^{y_0, y}$ then $T_0^X - T_0^Y$ converges in distribution to 0 as $y \uparrow y_0$. This and the continuity of $X$ under $\widetilde{\mathbb{Q}}_0^{y_0, y}$ imply that we can find $y_2 \in (0, y_0)$ close to $y_0$ so that for $y \in [y_2, y_0]$,

$$\mathbb{Q}_0^{y_0, y}(X_{T_0^Y} < y_1) \geq 1 - \frac{\delta y_0}{8}. \quad (3.35)$$

By Lemma 3.6 and scaling, there exists $y_3 \in (0, 1)$ small so that for $x \in (0, y_3)$ and $y \in [-x, 0]$,

$$\mathbb{P}^{x,y}(Y_{T_{y_0}^X} \in [y_2, y_0] | T_{y_0}^X < T_0^X) \geq 1 - \frac{\delta y_0}{8}.$$

A routine application of the theory of Doob's $h$-processes shows that the last estimate is equivalent to

$$\mathbb{Q}_0^{x,y}(Y_{T_{y_0}^X} \in [y_2, y_0] | T_{y_0}^X < T_0^X) \geq 1 - \frac{\delta y_0}{8}. \quad (3.36)$$

Let $S_1 = T_0^Y \circ T_{y_0}^X + T_{y_0}^X$ and $S_2 = T_0^X \circ S_1 + S_1$. We use the pseudo-strong Markov property and (3.31)–(3.36) to see that for $x \in (0, y_3)$ and $y \in [-x, 0]$,

$$\begin{aligned}
&\mathbb{Q}_0^{x,y}(Y_T \notin [-y_0, y_0]) \\
&= \mathbb{Q}_0^{x,y}(\{Y_T \notin [-y_0, y_0]\} \cap \{T_{y_0}^X < T_0^X\}) \\
&\leq \mathbb{Q}_0^{x,y}(T_{y_0}^X < T_0^X)\mathbb{Q}_0^{x,y}(Y_T \notin [-y_0, y_0] | T_{y_0}^X < T_0^X) \\
&\leq \frac{x}{y_0} \mathbb{Q}_0^{x,y}(Y_T \notin [-y_0, y_0] | T_{y_0}^X < T_0^X)
\end{aligned}$$

$$\leq \frac{x}{y_0} \mathbb{Q}_0^{x,y}(\{Y_{T_{y_0}^X} \notin [y_2, y_0]\} | T_{y_0}^X < T_0^X)$$
$$+ \frac{x}{y_0} \mathbb{Q}_0^{x,y}(\{Y_{T_{y_0}^X} \in [y_2, y_0]\} \cap \{X_{S_1} \geq y_1\} | T_{y_0}^X < T_0^X)$$
$$+ \frac{x}{y_0} \mathbb{Q}_0^{x,y}(\{Y_{T_{y_0}^X} \in [y_2, y_0]\} \cap \{X_{S_1} < y_1\} \cap \{Y_{S_2} \notin [-y_0, y_0]\} | T_{y_0}^X < T_0^X)$$
$$\leq \frac{x}{y_0} \mathbb{Q}_0^{x,y}(\{Y_{T_{y_0}^X} \notin [y_2, y_0]\} | T_{y_0}^X < T_0^X)$$
$$+ \frac{x}{y_0} \mathbb{Q}_0^{x,y}(\mathbf{1}_{\{Y_{T_{y_0}^X} \in [y_2, y_0]\}} \mathbb{Q}_0^{x,y}(X_{S_1} \geq y_1 | \mathcal{F}_{T_{y_0}^X}) | T_{y_0}^X < T_0^X)$$
$$+ \frac{x}{y_0} \mathbb{Q}_0^{x,y}(\mathbf{1}_{\{Y_{T_{y_0}^X} \in [y_2, y_0], X_{S_1} < y_1\}} \mathbb{Q}_0^{x,y}(Y_{S_2} \notin [-y_0, y_0] | \mathcal{F}_{S_1}) | T_{y_0}^X < T_0^X)$$
$$\leq \frac{x}{y_0} \left( \frac{\delta y_0}{8} + \frac{\delta y_0}{8} + \frac{\delta y_0}{8} \right)$$
$$\leq \delta x / 2.$$

As $|Y_T| \leq 1$ under $\mathbb{Q}_0^{x,y}$ by Lemma 3.8, we conclude from this and (3.30) that for $x \in (0, y_3)$ and $y \in [-x, 0]$,

$$(3.37) \quad \mathbb{E}^{x,y}[-Y_T \mathbf{1}_{\{|Y_T| \leq y_0\}} | X_T = 0] \geq \left(1 - \frac{\delta}{2}\right) x - \frac{\delta}{2} x = (1 - \delta) x.$$

Recall that in the above argument, we start with an arbitrary $y_0 \in (0, 1/2)$, then we find $y_1$, $y_2$ and $y_3$ accordingly so that (3.34), (3.35) and (3.36) hold, respectively. Let $r_0$ be the value of $y_3$ corresponding to $y_0 = 1/3$. We now define for $x \in (0, r_0]$,

$$\psi(x) = \inf \left\{ a > 0 : \inf_{y \in [-x, 0]} \mathbb{E}^{x,y}[-Y_T \mathbf{1}_{\{|Y_T| \leq a\}} | X_T = 0] \geq (1 - \delta) x \right\}.$$

Clearly by (3.37), $\psi(x) \leq 1/3$. Since for every $y_0 \in (0, 1/2)$, there is $y_3 > 0$ so that (3.37) holds, it follows that $\psi(x) \leq y_0$ for $x \in (0, y_3)$. Hence $\lim_{x \downarrow 0} \psi(x) = 0$. Summarizing, we obtain for $x \in (0, r_0)$ and $y \in [-x, 0]$,

$$\mathbb{E}^{x,y}[-Y_T \mathbf{1}_{\{|Y_T| \leq \psi(x)\}} | X_T = 0] \geq (1 - \delta) x.$$

This proves (3.28) and so (3.29) follows since, as we observed previously, $\mathbb{P}^{x,y}(X_T = 0) = 1 - x$. □

## 4. Pathwise uniqueness.

PROOF OF THEOREM 1.2. As we noticed previously, it suffices to prove pathwise uniqueness for solutions $X$ and $Y$ of (1.1) and (1.3) with $X_0 = Y_0 = 0$. The proof will be divided into three parts. The main argument will be presented in Part 1 and subdivided into three steps.

*Part* 1 (*Strong uniqueness*). We first show that strong uniqueness holds for solutions of (1.2)–(1.3) when there is a single filtration. Let $(X, W)$ and $(Y, W)$ be two weak solutions of (1.1) satisfying (1.3) with a common Brownian motion $W$ and such that $X_0 = Y_0 = 0$.

Define
$$M_t = |X_t| \vee |Y_t| \quad \text{and} \quad Z_t = |X_t - Y_t|.$$

Our strategy is to show that, with probability one, on any excursion of $M$ away from 0 that reaches level 1, $X$ and $Y$ have to agree. A scaling argument then shows that for any $b > 0$, with probability one, on any excursion of $M$ away from 0 that reaches level $b$, $X$ and $Y$ have to agree, and this will give the strong uniqueness for solutions of (1.2)–(1.3). We execute this plan in three steps.

*Step* 1. In this step, we show that if there exist two solutions, neither of them can stay on one side of 0 between the bifurcation time and the time when $M$ reaches level 1. In other words, setting

$$S = \inf\{t > 0 : M_t = 1\}, \qquad L = \sup\{t < S : M_t = 0\},$$
$$C_0 = \{X_s X_t > 0 \text{ for all } s, t \in (L, S]\} \cup \{Y_s Y_t > 0 \text{ for all } s, t \in (L, S]\},$$

we show that

(4.1) $$\mathbb{P}^{0,0}(\{\exists t \in [L, S] : X_t \neq Y_t\} \cap C_0) = 0.$$

For $b \in [0, 1)$ let
$$L_b = \inf\{t \geq L : M_t = b\},$$

and for $b \in (0, 1)$,

$$C_b := (\{|X_{L_b}| \geq |Y_{L_b}|\} \cap \{X_s X_t > 0 \text{ for } s, t \in (L_b, S]\})$$
$$\cup (\{|X_{L_b}| \leq |Y_{L_b}|\} \cap \{Y_s Y_t > 0 \text{ for } s, t \in (L_b, S]\}).$$

Note that $L_b$ is not a stopping time and thus we cannot apply the pseudo-strong Markov property at $T_b$ to estimate the probability of $C_b$. To circumvent this difficulty, we define two sequence of stopping times $\{T_j, j \geq 0\}$ and $\{S_j, j \geq 1\}$ as follows. Let $T_0 = 0$, and for $j \geq 1$,

$$S_j = \inf\{t > T_{j-1} : M_t = b\} \quad \text{and} \quad T_j = \inf\{t > S_j : M_t = 0\}.$$

It is possible that some or all of the above stopping times are infinite. For $a > 0$ set $U_a^{j,|X|} = \inf\{t > S_j : |X_t| = a\}$ and define $U_a^{j,|Y|}$ similarly. Let

$$A_{j,b} := \{T_{j-1} < S, S_j < \infty, |Y_{S_j}| \leq |X_{S_j}| = b \text{ and } U_1^{j,|X|} < U_0^{j,|X|}\}$$
$$\cup \{T_{j-1} < S, S_j < \infty, |X_{S_j}| \leq |Y_{S_j}| = b \text{ and } U_1^{j,|Y|} < U_0^{j,|Y|}\}.$$

On $A_{j,b}$, we have $L \in [T_{j-1}, S_j]$ and $L_b = S_j$. Thus the $\{A_{j,b}, j \geq 1\}$ are disjoint and $A_{j,b} \subset C_b$ for every $j \geq 1$. In particular we have $\bigcup_{j=1}^{\infty} A_{j,b} \subset C_b$. On the other hand, since $M$ is a continuous process, during any finite time interval, it can only oscillate between $0$ and $b$ a finite number of times. This implies that $C_b \subset \bigcup_{j=1}^{\infty} A_{j,b}$. Therefore we have for $0 < b < 1$,

$$(4.2) \qquad C_b = \bigcup_{j=1}^{\infty} A_{j,b}.$$

Applying the pseudo-strong Markov property at time $S_j$ and using Lemma 3.6 and symmetry, we can choose $b \in (0, 1)$ small enough so that for every $j \geq 1$,

$$\mathbb{P}^{0,0}(\{|Y_{U_1^{j,|X|}} - X_{U_1^{j,|X|}}| > \varepsilon\} \cap A_{j,b} \cap \{|Y_{S_j}| \leq |X_{S_j}| = b\})$$

$$= \mathbb{P}^{0,0}(\mathbb{P}^{0,0}(\{|Y_{U_1^{j,|X|}} - X_{U_1^{j,|X|}}| > \varepsilon\} \cap A_{j,b} \cap \{|Y_{S_j}| \leq |X_{S_j}| = b\}|\mathcal{F}_{S_j}))$$

$$= \mathbb{E}^{0,0}(\mathbf{1}_{\{T_{j-1}<S, S_j<\infty, |Y_{S_j}| \leq |X_{S_j}|=b\}} \mathbb{P}^{0,0}(|Y_{U_1^{j,|X|}} - X_{U_1^{j,|X|}}| > \varepsilon,$$

$$U_1^{j,|X|} < U_0^{j,|X|} | \mathcal{F}_{S_j}))$$

$$\leq \mathbb{E}^{0,0}(\mathbf{1}_{\{T_{j-1}<S, S_j<\infty, |Y_{S_j}| \leq |X_{S_j}|=b\}} \varepsilon \mathbb{P}^{0,0}(U_1^{j,|X|} < U_0^{j,|X|} | \mathcal{F}_{S_j}))$$

$$= \varepsilon \mathbb{P}^{0,0}(A_{j,b} \cap \{|Y_{S_j}| \leq |X_{S_j}| = b\}).$$

It follows that

$$\mathbb{P}^{0,0}(|Y_{U_1^{j,|X|}} - X_{U_1^{j,|X|}}| > \varepsilon | A_{j,b} \cap \{|Y_{S_j}| \leq |X_{S_j}| = b\}) \leq \varepsilon,$$

and, similarly,

$$\mathbb{P}^{0,0}(|X_{U_1^{j,|Y|}} - Y_{U_1^{j,|Y|}}| > \varepsilon | A_{j,b} \cap \{|X_{S_j}| \leq |Y_{S_j}| = b\}) \leq \varepsilon.$$

This implies that

$$\mathbb{P}^{0,0}(|X_S - Y_S| \leq \varepsilon | A_{j,b}) \geq 1 - \varepsilon$$

and so

$$\mathbb{P}^{0,0}(\{|X_S - Y_S| \leq \varepsilon\} \cap A_{j,b}) \geq (1 - \varepsilon) \mathbb{P}^{0,0}(A_{j,b}).$$

Summing over $j \geq 1$ yields

$$(4.3) \qquad \mathbb{P}^{0,0}(\{|X_S - Y_S| \leq \varepsilon\} \cap C_b) \geq (1 - \varepsilon) \mathbb{P}^{0,0}(C_b).$$

For $1 > b_1 > 4b_2 > 0$, in view of (4.2) for $C_{b_2}$ (with $b_2$ in place of $b$ there), we have by Corollary 3.9 and the pseudo-strong Markov property applied at stopping times $S_j$ that $C_{b_2} \subset C_{b_1}$, $\mathbb{P}^{0,0}$-a.s. Therefore $\bigcap_{n \geq 1} C_{5^{-n}} = C_0$, $\mathbb{P}^{0,0}$-a.s., and $\lim_{n \to \infty} \mathbb{P}^{0,0}(C_{5^{-n}}) = \mathbb{P}^{0,0}(C_0)$.

We can choose $b = 5^{-n}$ sufficiently small so that $\mathbb{P}^{0,0}(C_b) \geq (1-\varepsilon)\mathbb{P}^{0,0}(C_0)$. This and (4.3) imply that

$$\begin{aligned}
\mathbb{P}^{0,0}(\{|X_S - Y_S| \leq \varepsilon\} \cap C_0) \\
\geq \mathbb{P}^{0,0}(\{|X_S - Y_S| \leq \varepsilon\} \cap C_b) - \varepsilon\mathbb{P}^{0,0}(C_0) \\
\geq (1-\varepsilon)\mathbb{P}^{0,0}(C_b) - \varepsilon\mathbb{P}^{0,0}(C_0) \\
\geq ((1-\varepsilon)^2 - \varepsilon)\mathbb{P}^{0,0}(C_0).
\end{aligned}$$

Since $\varepsilon > 0$ is arbitrarily small, it follows that $X_S = Y_S$, $\mathbb{P}^{0,0}$-a.s., on $C_0$. It follows from Remark 2.3 that $X_t = Y_t$ for every $t \in [L, S]$, $\mathbb{P}^{0,0}$-a.s., on $C_0$. This proves (4.1).

For $b > 0$, let

$$S^b = \inf\{t > 0 : M_t = b\}, \qquad L^b = \sup\{t < S^b : M_t = 0\},$$
$$C_0^b = \{X_s X_t > 0 \text{ for } s, t \in (L^b, S^b]\} \cup \{Y_s Y_t > 0 \text{ for } s, t \in (L^b, S^b]\}.$$

Let $R_0^b = 0$, $\widehat{S}_1^b = S^b$, $\widehat{R}_1^b := \inf\{t > S^b : M_t = 0\}$, and for $k \geq 1$ define

$$\widehat{S}_k^b = \widehat{R}_{k-1}^b + S^b \circ \theta_{\widehat{R}_{k-1}^b},$$
$$\widehat{R}_k^b = \widehat{S}_{k-1}^b + R_1^b \circ \theta_{\widehat{S}_{k-1}^b},$$
$$\widehat{L}_k^b = \sup\{t < S_k^b : M_t = 0\},$$
$$\widehat{C}_{0,k}^b = \{X_s X_t > 0 \text{ for } s, t \in (L_k^b, S_k^b]\} \cup \{Y_s Y_t > 0 \text{ for } s, t \in (L_k^b, S_k^b]\}.$$

Then by the pseudo-strong Markov property applied at times $\widehat{R}_k^b$,

$$\mathbb{P}^{0,0}\left(\bigcup_{k \geq 1}\{\exists t \in [L_k^1, S_k^1] : X_t \neq Y_t\} \cap \widehat{C}_{0,k}^1\right) = 0.$$

In an analogous way we can prove that for any $b > 0$,

(4.4) $$\mathbb{P}^{0,0}\left(\bigcup_{k \geq 1}\{\exists t \in [L_k^b, S_k^b] : X_t \neq Y_t\} \cap \widehat{C}_{0,k}^b\right) = 0.$$

*Step* 2. In this intermediate step, we show that for any two arbitrary small constants $b, \varepsilon_0 > 0$, there is some $a_1 = a_1(b, \varepsilon) > 0$ such that for every $b_0 \in (0, a_1]$,

(4.5) $$\mathbb{P}^{b_0, 0}(T_b^Z < T_1^M | T_1^M \leq T_0^Z) \leq \varepsilon_0.$$

Fix arbitrarily small $b, \varepsilon_0 > 0$ and a large enough integer $m \geq 2$ such that $\frac{1}{(m-1)b} \leq \varepsilon_0$. Fix $\delta \in (0, 1/4)$ small such that

$$\sum_{j=0}^{m-1}(1-\delta)^{2j} \geq m-1.$$

So for every $a > 0$,

$$（4.6) \quad 1 - a \sum_{j=0}^{m-1} (1-\delta)^j (1-\delta)^j \leq 1 - (m-1)a.$$

Choose a constant $r_0 \in (0, \delta)$ and a function $\psi$ that satisfies the statement of Lemma 3.10 together with the given $\delta$. Make $r_0 > 0$ smaller, if necessary, so that $r_0 < 1/(2m)$. Let $a_1 \in (0, b/2)$ be small enough so that

$$\psi^m(a_1) := \underbrace{\psi \circ \psi \circ \cdots \circ \psi}_{m \text{ times}}(a_1) < r_0.$$

Assume that $X_0 = b_0 \in (0, a_1)$ and $Y_0 = 0$. Let $U_0 = 0$ and for $k \geq 1$,

$$U_k = \begin{cases} \inf\{t \geq U_{k-1} : X_t = 0\}, & \text{if } Y_{U_{k-1}} = 0, \\ \inf\{t \geq U_{k-1} : Y_t = 0\}, & \text{if } X_{U_{k-1}} = 0. \end{cases}$$

It follows from (4.2) that for $0 \leq n \leq m$, on the event $\{T_1^M > U_n\} \cap \{Z_{U_n} \leq \psi^n(b_0)\}$,

$$(4.7) \quad \mathbb{E}^{b_0,0}[Z_{U_{n+1}} \mathbf{1}_{\{Z_{U_{n+1}} \leq \psi^{n+1}(b_0)\}} \mathbf{1}_{\{T_1^M > U_{n+1}\}} | \mathcal{F}_{U_n}]$$

$$\geq (1-\delta) Z_{U_n}(1 - Z_{U_n}) \mathbf{1}_{\{Z_{U_n} \leq \psi^n(b_0)\}}$$

Since $X$ is a continuous local martingale, by the gambler's ruin estimate, we have for $n \geq 0$ that on the event $\{T_1^M > U_n\}$,

$$(4.8) \quad \mathbb{P}^{b_0,0}(T_1^M > U_{n+1} | \mathcal{F}_{U_n}) = 1 - Z_{U_n}.$$

Recall that $\psi^m(b_0) \leq r_0 < 1/(2m)$. Then, for $\gamma \in (0, m)$,

$$1 - \gamma(1-\delta) Z_{U_n} \mathbf{1}_{\{Z_{U_n} \leq \psi^n(b_0)\}} \geq 0.$$

These remarks, (4.7) and (4.8) imply that for $\gamma \in (0, m)$ and $n \leq m$,

$$\mathbb{E}^{b_0,0}[(1 - \gamma Z_{U_{n+1}} \mathbf{1}_{\{Z_{U_{n+1}} \leq \psi^{n+1}(b_0)\}}) \mathbf{1}_{\{T_1^M > U_{n+1}\}}]$$

$$= \mathbb{E}^{b_0,0}[(1 - \gamma Z_{U_{n+1}} \mathbf{1}_{\{Z_{U_{n+1}} \leq \psi^{n+1}(b_0)\}}) \mathbf{1}_{\{T_1^M > U_{n+1}\}} \mathbf{1}_{\{T_1^M > U_n\}}]$$

$$= \mathbb{E}^{b_0,0}[\mathbb{E}^{b_0,0}[(1 - \gamma Z_{U_{n+1}} \mathbf{1}_{\{Z_{U_{n+1}} \leq \psi^{n+1}(b_0)\}}) \mathbf{1}_{\{T_1^M > U_{n+1}\}} | \mathcal{F}_{U_n}] \mathbf{1}_{\{T_1^M > U_n\}}]$$

$$\leq \mathbb{E}^{b_0,0}[(1 - Z_{U_n}) \mathbf{1}_{\{T_1^M > U_n\}}]$$

$$\quad - \gamma \mathbb{E}^{b_0,0}[\mathbb{E}^{b_0,0}[Z_{U_{n+1}} \mathbf{1}_{\{Z_{U_{n+1}} \leq \psi^{n+1}(b_0)\}} \mathbf{1}_{\{T_1^M > U_{n+1}\}} | \mathcal{F}_{U_n}]$$

$$\quad\quad\quad\quad \times \mathbf{1}_{\{T_1^M > U_n\}} \mathbf{1}_{\{Z_{U_n} \leq \psi^n(b_0)\}}]$$

$$\leq \mathbb{E}^{b_0,0}[(1 - Z_{U_n} \mathbf{1}_{\{Z_{U_n} \leq \psi^n(b_0)\}}) \mathbf{1}_{\{T_1^M > U_n\}}]$$

$$\quad - \gamma \mathbb{E}^{b_0,0}[(1-\delta) Z_{U_n}(1 - Z_{U_n}) \mathbf{1}_{\{T_1^M > U_n\}} \mathbf{1}_{\{Z_{U_n} \leq \psi^n(b_0)\}}]$$

$$\leq \mathbb{E}^{b_0,0}[((1 - Z_{U_n}\mathbf{1}_{\{Z_{U_n} \leq \psi^n(b_0)\}}) - \gamma(1-\delta)(1-r_0)Z_{U_n}\mathbf{1}_{\{Z_{U_n} \leq \psi^n(b_0)\}})$$
$$\times \mathbf{1}_{\{T_1^M > U_n\}}]$$
$$= \mathbb{E}^{b_0,0}[(1 - (\gamma(1-\delta)(1-r_0)+1)Z_{U_n}\mathbf{1}_{\{Z_{U_n} \leq \psi^n(b_0)\}})\mathbf{1}_{\{T_1^M > U_n\}}].$$

An induction argument based on the above inequality shows that for $n \leq m$,

$$\mathbb{E}^{b_0,0}[(1 - Z_{U_{n-1}}\mathbf{1}_{\{Z_{U_{n-1}} \leq \psi^{n-1}(b_0)\}})\mathbf{1}_{\{T_1^M > U_{n-1}\}}]$$
$$\leq q\mathbb{E}^{b_0,0}\left[\left(1 - \left(\sum_{j=0}^{n-1}(1-\delta)^j(1-r_0)^j\right)Z_{U_0}\mathbf{1}_{\{Z_{U_0} \leq b_0\}}\right)\mathbf{1}_{\{T_1^M > U_0\}}\right]$$
$$= 1 - b_0 \sum_{j=0}^{n-1}(1-\delta)^j(1-r_0)^j$$
$$\leq 1 - (m-1)b_0.$$

We obtain from the above, (4.6) and (4.8),

$$\mathbb{P}^{b_0,0}(T_1^M > U_m) = \mathbb{E}^{b_0,0}[\mathbb{P}^{b_0,0}(T_1^M > U_m | \mathcal{F}_{U_{m-1}})\mathbf{1}_{\{T_1^M > U_{m-1}\}}]$$
$$= \mathbb{E}^{b_0,0}[(1 - Z_{U_{m-1}})\mathbf{1}_{\{T_1^M > U_{m-1}\}}]$$
$$\leq \mathbb{E}^{b_0,0}[(1 - Z_{U_{m-1}}\mathbf{1}_{\{Z_{U_{m-1}} \leq \psi^{n-1}(b_0)\}})\mathbf{1}_{\{T_1^M > U_{m-1}\}}]$$
$$\leq 1 - (m-1)b_0.$$

Thus

$$\mathbb{P}^{b_0,0}(T_1^M \leq T_0^Z) \geq \mathbb{P}^{b_0,0}(T_1^M \leq U_m) \geq (m-1)b_0.$$

Since $Z_{t \wedge T_0^Z}$ is a continuous local martingale, the gambler's ruin estimate tells us that

$$\mathbb{P}^{b_0,0}(T_b^Z \leq T_0^Z) \leq b_0/b,$$

and, therefore,

$$\mathbb{P}^{b_0,0}(T_b^Z < T_1^M | T_1^M \leq T_0^Z) = \frac{\mathbb{P}^{b_0,0}(T_b^Z < T_1^M \leq T_0^Z)}{\mathbb{P}^{b_0,0}(T_1^M \leq T_0^Z)} \leq \frac{\mathbb{P}^{b_0,0}(T_b^Z < T_0^Z)}{\mathbb{P}^{b_0,0}(T_1^M \leq T_0^Z)}$$
$$\leq \frac{1}{(m-1)b} \leq \varepsilon_0.$$

This proves (4.5). By symmetry, inequalities analogous to (4.5) hold when $(b_0, 0)$ is replaced by $(-b_0, 0)$, $(0, b_0)$ or $(0, -b_0)$.

*Step* 3. We complete the proof of the claim that with probability one, on any excursion of $M_t = |X_t| \vee |Y_t|$ away from 0 that reaches level 1, $Z_t =$ $|X_t - Y_t|$ must be zero. That is, using the definitions of $S$ and $L$ from Step 1, we will show in this step that

(4.9) $$\mathbb{P}^{0,0}(X_t = Y_t \text{ for every } t \in [L, S]) = 1.$$

Note that $S < \infty$ with probability one.

We chose $b, \varepsilon_0 > 0$ and $a_1 = a_1(b, \varepsilon_0) > 0$ in Step 2. Define for $0 < b_1 < a_1$,

$$F_{a_1,b_1} = \{\exists t \in (L, S) : Z_t = M_t \in (b_1, a_1]\},$$
$$F_{0+} = \{\forall u \in (L, S) \ \exists t \in (L, u) : Z_t = M_t\}.$$

It follows from (4.4) applied to all rational $b > 0$ that $\mathbb{P}^{0,0}$-a.s., for any (random) $0 \leq t_0 < t_1$ such that $M_{t_0} = Z_{t_0} = 0$ and $Z_t > 0$ for $t \in (t_0, t_1)$, there is a (random) decreasing sequence $\{t_n, n \geq 1\} \subset (t_0, t_1)$ such that $t_n \downarrow t_0$ as $n \to \infty$ and $M_{t_n} = Z_{t_n}$ for every $n \geq 1$. Thus to prove (4.9), it will suffice to show that

(4.10) $$\mathbb{P}^{0,0}(\{\exists t \in [L, S] : X_t \neq Y_t\} \cap F_{0+}) = 0.$$

Let $\tau_0 = 0$ and for $k \geq 1$, define

$$\sigma_k = \inf\{t > \tau_{k-1} : Z_t = M_t \in (b_1, a_1]\} \quad \text{and} \quad \tau_k = \inf\{t > \sigma_k : M_t = 0\}.$$

We further define

$$\tau_1^{M,k} = \inf\{t > \sigma_k : M_t \geq 1\} \quad \text{and} \quad \tau_b^{Z,k} = \inf\{t > \sigma_k : Z_t \geq b\}.$$

Note that

$$F_{a_1,b_1} = \bigcup_{k=1}^{\infty} \{\sigma_k < S \text{ and } \tau_1^{M,k} < \tau_k\}$$

$$= \bigcup_{k=1}^{\infty} \left( \bigcap_{j=1}^{k-1} \{\sigma_j < S \text{ and } \tau_1^{M,j} \geq \tau_j\} \cap \{\sigma_k < S \text{ and } \tau_1^{M,k} < \tau_k\} \right),$$

with the convention that $\bigcap_{j=1}^{0} \{\sigma_j < S \text{ and } \tau_1^{M,j} \geq \tau_j\} = \Omega$.

For every $k \geq 1$, applying the pseudo-strong Markov property at $\sigma_k$ and using Step 2, we have

$$\mathbb{P}^{0,0}\left( \tau_b^{Z,k} < \tau_1^{M,k} \,\Big|\, \bigcap_{j=1}^{k-1} \{\sigma_j < S \text{ and } \tau_1^{M,j} \geq \tau_j\} \cap \{\sigma_k < S \text{ and } \tau_1^{M,k} < \tau_k\} \right) \leq \varepsilon_0.$$

This implies that

$$\mathbb{P}^{0,0}\left( \{|Y_S - X_S| < b\} \cap \bigcap_{j=1}^{k-1} \{\sigma_j < S \text{ and } \tau_1^{M,j} \geq \tau_j\} \cap \{\sigma_k < S \text{ and } \tau_1^{M,k} < \tau_k\} \right)$$

$$\geq (1 - \varepsilon_0) \mathbb{P}^{0,0}\left( \bigcap_{j=1}^{k-1} \{\sigma_j < S \text{ and } \tau_1^{M,j} \geq \tau_j\} \cap \{\sigma_k < S \text{ and } \tau_1^{M,k} < \tau_k\} \right).$$

Summing over $k \geq 1$, we have

(4.11) $\quad \mathbb{P}^{0,0}(\{|Y_S - X_S| < b\} \cap F_{a_1,b_1}) \geq (1 - \varepsilon_0)\mathbb{P}^{0,0}(F_{a_1,b_1}).$

Note that $F_{a_1,b_1} \subset F_{a_1,b_2}$ if $b_1 > b_2$. Hence for sufficiently small $b_1 > 0$, $\mathbb{P}^{0,0}(F_{a_1,b_1}) \geq (1 - \varepsilon_0)\mathbb{P}^{0,0}(F_{a_1,0})$. This and (4.11) imply that

$$\begin{aligned}
\mathbb{P}^{0,0}&(\{|X_S - Y_S| \leq b\} \cap F_{a_1,0}) \\
&\geq \mathbb{P}^{0,0}(\{|X_S - Y_S| \leq b\} \cap F_{a_1,b_1}) - \varepsilon_0 \mathbb{P}^{0,0}(F_{a_1,0}) \\
&\geq (1 - \varepsilon)\mathbb{P}^{0,0}(F_{a_1,b_1}) - \varepsilon_0 \mathbb{P}^{0,0}(F_{a_1,0}) \\
&\geq ((1 - \varepsilon_0)^2 - \varepsilon_0)\mathbb{P}^{0,0}(F_{a_1,0}).
\end{aligned}$$
(4.12)

Note that $F_{a_1,0} \subset F_{a_2,0}$ if $a_1 < a_2$. If the event $\bigcap_{a>0} F_{a,0}$ holds, then there exist $t_n \in (L, S)$ such that $Z_{t_n} = M_{t_n} \in (0, 1/n]$ for all $n \geq 1$. By compactness, $t_n$ must have a subsequence $t_{n_k}$ converging to a point $t_\infty \in [L, S]$. By the continuity of $X$ and $Y$, $Z_{t_\infty} = M_{t_\infty} = 0$, so it follows from the definition of $L$ that $t_\infty = L$. We conclude that $\bigcap_{a>0} F_{a,0} = F_{0+}$, $\mathbb{P}^{0,0}$-a.s. Thus, for sufficiently small $a_1 > 0$, $\mathbb{P}^{0,0}(F_{0+}) \geq (1 - \varepsilon_0)\mathbb{P}^{0,0}(F_{a_1,0})$. This and (4.12) imply that for every $\varepsilon_0 > 0$,

$$\begin{aligned}
\mathbb{P}^{0,0}&(\{|X_S - Y_S| \leq b\} \cap F_{0+}) \\
&\geq \mathbb{P}^{0,0}(\{|X_S - Y_S| \leq b\} \cap F_{a_1,0}) - \varepsilon_0 \mathbb{P}^{0,0}(F_{a_1,0}) \\
&\geq \mathbb{P}^{0,0}(\{|X_S - Y_S| \leq b\} \cap F_{a_1,0}) - (\varepsilon_0/(1 - \varepsilon_0))\mathbb{P}^{0,0}(F_{0+}) \\
&\geq ((1 - \varepsilon_0)^2 - \varepsilon_0)\mathbb{P}^{0,0}(F_{a_1,0}) - (\varepsilon_0/(1 - \varepsilon_0))\mathbb{P}^{0,0}(F_{0+}) \\
&\geq ((1 - \varepsilon_0)^2 - \varepsilon_0 - (\varepsilon_0/(1 - \varepsilon_0)))\mathbb{P}^{0,0}(F_{0+}).
\end{aligned}$$

Since $\varepsilon_0 > 0$ and $b > 0$ are arbitrarily small, it follows that $\mathbb{P}^{0,0}(\{X_S \neq Y_S\} \cap F_{0+}) = 0$. In view of Remark 2.3 and Step 1, this proves (4.9). Another application of Remark 2.3 and (4.9) yields

(4.13) $\quad \mathbb{P}^{0,0}(X_t = Y_t \text{ for every } t \in [L, R]) = 1,$

where $R = \inf\{t > S : M_t = 0\}$.

Recall the definitions of $S_k^b, R_k^b$ and $L_k^b$ from Step 1. Just as in the case of (4.4), we can deduce from (4.13) that for $b > 0$,

$$\mathbb{P}^{0,0}\left(\bigcup_{k \geq 1}\{\exists t \in [L_k^b, S_k^b] : X_t \neq Y_t\}\right) = 0.$$

Since $b$ can be arbitrarily small, this proves that $\mathbb{P}^{0,0}$-a.s., $X_t = Y_t$ for every $t \geq 0$.

So far, our entire proof was concerned with processes defined on the canonical space. Now suppose that $X'$ and $Y'$ are two weak solutions to (1.1) and (1.3) driven by the same Brownian motion $W'$, starting from $W'_0 = 0$, $X'_0 = 0$ and $Y'_0 = 0$, and defined on some probability space $(\Omega', \mathbb{P}')$. Using the definition of $\mathbb{P}^{0,0}$ given in (2.1) it is clear that $X'_t = Y'_t$ for every $t \geq 0$, $\mathbb{P}'$-a.s.

*Part* 2 (*Strong existence*). Existence of a weak solution to the SDE (1.2) that spends zero time at 0 follows from [7]. The solution can be constructed as a time change of Brownian motion. Existence of a strong solution for the SDE (1.2) that spends zero time at 0 follows from the strong uniqueness and weak existence for solutions of (1.2) that spend zero time at 0; this can be done in the same way as in [12], or following word-for-word the proof of [11], Theorem IX.1.7(ii).

*Part* 3 (*Pathwise uniqueness*). Note that strong uniqueness implies weak uniqueness by the proof in [11], Theorem IX.1.7(i). Pathwise uniqueness now follows from strong existence and strong uniqueness by the same argument as in the last paragraph of the proof of [4], Theorem 5.8, or in the last paragraph of the proof of [5], Theorem 5.3. □

## 5. Stochastic differential equations with reflection.

PROOF OF THEOREM 1.3. Note that under the assumptions of Theorem 1.3, $a(x)^{-2}$ and hence $b(x)a(x)^{-2}$ are locally integrable on $\mathbb{R}$, and, therefore, $\int_0^y b(r)a(r)^{-2}\,dr$ is a continuous strictly increasing function. Define

$$s(x) := \int_0^x \exp\left(-\int_0^y \frac{2b(r)}{a(r)^2}\,dr\right)dy, \qquad x \in \mathbb{R}.$$

We will use $s^{-1}$ to denote the inverse function of $s$. If $(X, W)$ is a weak solution to (1.4)–(1.5), by the Itô–Tanaka formula (see [11], Theorem VI.1.5), we have

$$ds(X_t) = s'(X_t)a(X_t)\,dW_t + s'(X_t)\,dL_t.$$

Then $U := s(X)$ spends zero time at 0 and solves

(5.1) $\quad dU_t = (s'a) \circ s^{-1}(U_t)\,dW_t + s' \circ s^{-1}(U_t)\,dL_t \qquad$ with $U_0 = s(X_0)$.

By the uniqueness of the deterministic Skorokhod problem on $[0, \infty)$, there is a weak solution $U$ to (5.1) that spends zero time at 0, obtained as a time change of reflecting Brownian motion on $[0, \infty)$; moreover weak uniqueness holds for solutions of (5.1) that spend zero time at 0 (cf. [5], Section 4). It follows then that weak existence and weak uniqueness holds for solutions of (1.4) and (1.5).

Let $X$ and $Y$ be two weak solutions to (1.4)–(1.5) with the same driving Brownian motion $W$ with respect to a common filtration on a common probability space. Using the Itô–Tanaka formula,

$$(X_t - Y_t)^+ = (X_0 - Y_0)^+ + \int_0^t \mathbf{1}_{\{X_s - Y_s > 0\}}\,d(X_s - Y_s) + M_t,$$

where $M_t$ is a continuous nondecreasing process that increases only when $X_s - Y_s = 0$. Since $X_t \vee Y_t = Y_t + (X_t - Y_t)^+$, then

$$X_t \vee Y_t = X_0 \vee Y_0 + \int_0^t a(X_s \vee Y_s)\,dW_s + \int_0^t b(X_s \vee Y_s)\,ds + A_t,$$

where

$$A_t = \int_0^t \mathbf{1}_{\{X_s - Y_s > 0\}}\,d(L_s^X - L_s^Y) + L_t^Y + M_t.$$

When $X_s > Y_s$, then $dA_s = dL_s^X + dM_s$. This is 0 because $X_s \neq Y_s$, so that $dM_s = 0$, and $X_s > Y_s \geq 0$, so that $dL_s^X = 0$. When $X_s < Y_s$, $dA_s = dL_s^Y + dM_s$, which is 0 for the same reasons. When $X_s = Y_s = 0$, then $dA_s = dL_s^Y + dM_s \geq 0$. Finally, when $X_s = Y_s > 0$, then $dA_s = dM_s$. However, the argument of Le Gall [8] shows that the local time at 0 of $(X_t - Y_t)^+$ is 0 when $X_t$ and $Y_t$ are both in an interval for which either condition (a) or (b) holds. By our assumptions on $a$, this will be true when $X_t$ and $Y_t$ are both in any closed interval not containing 0. Therefore $A_t$ is nondecreasing and increases only when $X_t = Y_t = 0$.

If we let $Z_t = X_t \vee Y_t$, we then see that $Z_t$ is again a weak solution to (1.4) driven by the Brownian motion $W$ that spends zero time at 0. By the weak uniqueness for solutions of (1.4) that spend zero time at 0, the law of $Z_t$ is the same as that of $X_t$ and $Y_t$. But $Z_t \geq X_t$ for all $t$. We conclude that $Z_t = X_t$ for all $t$, and the same is true with $X$ replaced by $Y$. Therefore $X_t = Y_t$ for every $t \geq 0$. This proves the strong uniqueness for solutions of (1.4) that spend zero time at 0. The existence of a strong solution for SDE (1.4) that spends zero time at 0 follows from the strong uniqueness and weak existence for solutions of (1.4) that spend zero time at 0 in the same way as in [12], or as in the proof for [11], Theorem IX.1.7(ii). Pathwise uniqueness then follows in the same way as in the last paragraph of the proof of Theorem 1.2. □

PROOF OF COROLLARY 1.4. Suppose that functions $a$ and $b$ satisfy the assumptions of the corollary. Let $(X, W)$ and $(\widetilde{X}, \widetilde{W})$ be two weak solutions to (1.8)–(1.9) with a common Brownian motion $W$ (relative to possibly different filtrations) on a common probability space and $X_0 = \widetilde{X}_0$. As $a$ is an odd function and $X$ satisfies (1.9), by Tanaka's formula,

(5.2) $$d|X_t| = a(|X_t|)\,dW_t + b(|X_t|)\,dt + L_t,$$

where $L_t$ is the symmetric local time of $X$ which increases only when $X_t = 0$. Similarly, $|\widetilde{X}|$ satisfies equation (5.2) with $|\widetilde{X}|$ in place of $|X|$. Applying Theorem 1.3 to $|X|$ and $|\widetilde{X}|$, we have

$$\mathbb{P}(|X_t| = |\widetilde{X}_t| \text{ for all } t \geq 0) = 1. \qquad \Box$$

REMARK 5.1. (1) One reason the proof of Theorem 1.3 is considerably easier than that of Theorem 1.2 is that any two candidate solutions must be on the same side of 0. We tried to apply the method of proof of Theorem 1.3 to Theorem 1.2, but were unsuccessful.

(2) The function $a(x) = 2 + \sin(1/x^4)$ is an example of a function satisfying the hypotheses of Theorem 1.3. Because $a$ is bounded below away from 0, it is easy to show that any solution to the stochastic differential equation will spend zero time at 0. If we replace $|x|^\alpha$ in (1.2) by this $a$, will there be pathwise uniqueness in the two-sided case?

**Acknowledgments.** We are grateful to Tokuzo Shiga for bringing his joint paper [9] with S. Manabe to our attention and to Yuta Taniguchi for pointing out an error in a previous version of this paper.

R. F. Bass  
Department of Mathematics  
University of Connecticut  
Storrs, Connecticut 06269-3009  
USA  
E-mail: bass@math.uconn.edu

K. Burdzy  
Z.-Q. Chen  
Department of Mathematics  
University of Washington  
Seattle, Washington 98195  
USA  
E-mail: burdzy@math.washington.edu  
        zchen@math.washington.edu